\documentclass[a4paper,11pt]{amsart}
\usepackage[latin1]{inputenc}

\usepackage{multicol}
\usepackage{amssymb, amsmath, amsthm}
\usepackage{latexsym}
\usepackage{dsfont}
\usepackage{graphicx}

\usepackage{curves} 
\usepackage{epic} 
\usepackage{eepic}
\usepackage{xcolor}
 
\input xy
\xyoption{all}

\usepackage{verbatim}

\usepackage{hyperref}

\newtheorem{theorem}{Theorem}[section]

\newtheorem{lemma}[theorem]{Lemma}

\newtheorem{proposition}[theorem]{Proposition}

	{\par\noindent{\bf Proposition \ref{res:hiper}.}\!\!
	\nopagebreak\normalsize\it}%

	{\par\noindent{\bf Theorem \ref{result43}.}\!\!
	\nopagebreak\normalsize\it}%

	{\par\noindent{\it Sketch of the proof}.  
	\nopagebreak\normalsize}%
	{\hfill\linebreak[2]\hspace*{\fill}$\circlearrowleft$}

	{\par\noindent{\it Proof of Proposition }\ref{prop:stab:smc}.  
	\nopagebreak\normalsize}%
	{\hfill\linebreak[2]\hspace*{\fill}$\circlearrowleft$}

	{\par\noindent{\it Proof of Propositions }\ref{adap:mon}{\it and }\ref{simult:adap}.\!\!\!
	\nopagebreak\normalsize}%
	{\hfill\linebreak[2]\hspace*{\fill}$\circlearrowleft$}

{\theoremstyle{definition}
       \newtheorem{definition}[theorem]{Definition}

       \newtheorem{remark}[theorem]{Remark}
       
       \newtheorem{parrafo}[theorem]{{\!}}  }

\numberwithin{equation}{theorem}

\newcommand{\nat}{\mathbb N}
 
\newcommand{\calo}{{\mathcal {O}}}

\DeclareMathOperator{\ord}{ord}

\DeclareMathOperator{\Sing}{Sing}
\DeclareMathOperator{\Spec}{Spec}

\DeclareMathOperator{\Ord}{Ord}
\DeclareMathOperator{\Hord}{H-ord}
\DeclareMathOperator{\In}{In}
\DeclareMathOperator{\Gr}{Gr}
\DeclareMathOperator{\inv}{inv}

\newcommand{\R}{{\mathcal R}}

\newcommand{\G}{{\mathcal G}}
\newcommand{\C}{{\mathcal C}}

\renewcommand{\L}{{\mathcal L}}

\newcommand{\N}{{\mathbb N}}

\newcommand{\m}{{\mathcal{M}}}
\renewcommand{\P}{\mathcal{P}}
\newcommand{\Q}{\mathcal{Q}}
\newcommand{\A}{\mathbb{A}}

\newcommand{\id}[1]{\langle #1 \rangle}
\newcommand{\x}{{\mathbf{x}}}

\definecolor{darkpurple}{rgb}{0.28,0.24,0.55}
\definecolor{lightblue}{rgb}{0,0.75,1}

\title[Techniques for the study of singularities...]{Techniques for the study of singularities with applications to resolution of $2$-dimensional schemes}

\author{Ang\'elica Benito}
\author{Orlando E. Villamayor U.}
 \thanks{The authors are partially supported by MTM2009-07291.}
\date{\today}

\address{Dpto. Matem\'aticas,  Universidad
Aut\'onoma de Madrid and Instituto de Ciencias Matem\'aticas CSIC-UAM-UC3M-UCM \\
Ciudad Universitaria de Cantoblanco, 28049 Madrid, Spain}

\email[Ang\'elica Benito]{angelica.benito@uam.es}
\email[Orlando E. Villamayor U.]{villamayor@uam.es}


\begin{document}

\maketitle

\vspace{-.6cm}

\begin{abstract} 
We give an overview of invariants of algebraic singularities over perfect fields. We then show how they lead to a synthetic proof of embedded resolution of singularities of $2$-dimensional schemes.
\end{abstract}

\vspace{-.2cm}
{\tableofcontents}

\section{Introduction}
\begin{parrafo}
This paper includes an exposition of recent progress concerning singularities over perfect fields, it is also shown how these results lead to the resolution of singularities of $2$-dimensional schemes. 

The ultimate motivation of the results reported in this work is the open problem of resolution of singularities in any dimension. We focus here on resolution in the sense of Hironaka, which is a step by step procedure: Namely, given a reduced scheme $X$ over a perfect field $k$, the question is to construct a sequence of blow-ups along smooth centers, each center included in the Hilbert-Samuel stratum of the successive strict transforms of $X$, so as to define a desingularization. 

Sections \ref{sec2} and  \ref{sec3} include an overview of invariants that have been introduced recently for singularities over perfect fields in arbitrary dimension, and finally in Section \ref{sec4} we apply them to give a synthetic proof of resolution of $2$-dimensional schemes.

In the case of characteristic zero, given $X$ included in a smooth scheme of dimension $d$, say $V^{(d)}$, the existence of smooth hypersurfaces of {\em maximal contact} was used by Hironaka to attach to $X$ inductive invariants in dimension $d-1$. This is done  in  a way that the problem of resolution, formulated in dimension $d$, is reformulated as a problem in dimension $d-1$ by {\em restriction} to smooth hypersurfaces of maximal contact. A drawback of this approach is that such nice hypersurfaces are not unique, and, in addition to that, these hypersurfaces are defined only locally whereas the problem of resolution is global. So a significant  difficulty in using this form of induction is that of {\em patching} local information. 

This significant difficulty does not show up in Hironaka's theorem in \cite{Hir64}, which is existential, but it does arise in the {\em constructive} proofs of Hironaka's theorem. Namely in proofs which establish an algorithm that indicate which is the smooth center to be blown up in Hironaka's step by step procedure of resolution (see \cite{Villa89}, \cite{Villa92}, \cite{BM97}). 

A second drawback in using hypersurfaces of maximal contact for inductive arguments appears when trying to work over fields of positive characteristic. In fact, in this context these hypersurfaces do not always exist (\cite{Nar2}).

Section \ref{sec2} is an exposition of an approach to resolution, developed in recent years, which avoids the use of maximal contact. In fact, smooth hypersurfaces  of maximal contact can be replaced by local projections on smooth schemes; and 
the traditional notion of {\em restriction to hypersurfaces of maximal contact} can be replaced by a generalized form of the discriminant, defined in terms of elimination theory. 
When considering resolution of singularities over fields of characteristic zero this alternative procedure provides the same information as that obtained by using maximal contact. But it has an important advantage over the original constructive proofs of resolution as it trivializes the globalization of local invariants. This clarifies, in particular, the globalization of  local data extracted from the Hilbert Samuel function.  All this has led to a significant conceptual simplification of constructive resolution of singularities in characteristic zero (\cite{BV3}).

%
%
%

A second advantage of this approach is the fact that projections are definable over perfect fields; a feature which has opened the way to the definition of
new inductive invariants in positive characteristic. 
This study is addressed in Section \ref{sec3}, which is devoted to the discussion of inductive invariants of singularities over perfect fields, introduced in \cite{BeV1} and \cite{BeV2}. These are  natural extensions of the inductive invariants used by Hironaka in characteristic zero.

As was already mentioned, Hironaka uses inductive invariants in dimension $d-1$, in his step by step procedure, to obtain a resolution of singularities by successive monoidal transformations. These invariants enable him to construct a sequence of blow-ups over $X\subset V^{(d)}$, so as to come to a so called ``$d-1$-simplification''. This $d-1$-simplification is also known as a reduction to the {\em monomial case}. He then shows that 
it is easy to achieve desingularization once $X$ has been transformed into a scheme with singularities in the monomial case.

The inductive invariants, used by Hironaka (in dimension $d-1$) make use of hypersurfaces of maximal contact, and the argument works exclusively in characteristic zero. The alternative approach to induction, using projections, has led to the construction of a sequence of blow-ups over $X\subset V^{(d)}$, so as to come to a ``$d-1$-simplification'' now over perfect fields. This parallels the $d-1$-simplification in characteristic zero, or say the \emph{reduction to the monomial case}, but the outcome obtained is weaker in positive characteristic. 
In fact, over fields of positive characteristic the resolution of singularities which are in the monomial case is not straightforward, as it is in the case of characteristic zero. Despite this fact, the reduction to the monomial case is expected to be a simplification of the singularities. The main outcome of Section \ref{sec3} is to show that, if some  additional numerical conditions are fulfilled, singularities which are in the monomial case can be resolved. 

As an application, in Section \ref{sec4}, we show how this alternative approach to induction leads to resolution in the case in which $X\ (\subset V^{(d)})$ is a $2$-dimensional embedded scheme. 
The task, for future research, would be to show that the said numerical conditions, discussed in Section \ref{sec3}, can be attained in higher dimensions.

The reduction to the monomial case in positive characteristic is related with other forms of simplification that appear in other works. For example,  a procedure is introduced in the work of Kawanoue and Matsuki (\cite{Kaw} and \cite{KM}), which also parallels this reduction. See also \cite{Hironaka06} and \cite{Wlo2}.

The general and unifying strategy along this paper is the use of higher differential operators as a tool for the study of singularities 
over perfect field. Embedded resolution leads to the study of ideals on smooth schemes, and more generally to Rees algebras 
of ideals on smooth schemes. 
This reformulation of resolution problems appears already in Hironaka's work. Rees algebras on smooth schemes can be enriched by the action of higher order differential operators, and these enriched algebras are called differential algebras. These are algebras which encode very subtle information of the singularities, and they are the main tool in this alternative form of induction  which avoids the use of maximal contact. 

Differential Rees algebras allow us to reduce the problem of the resolution of singularities embedded in a smooth $d$-dimensional scheme, say $V^{(d)}$, to that of hypersurfaces of multiplicity 
$p^e$, where $p$ denotes the characteristic of the underlying field. More precisely to hypersurfaces defined by equations of the form
$$f_{p^{e}}(z_1)=z^{p^{e}}+a_1^{}z_1^{p^{e}-1}+\dots+a_{p^{e}}^{}\in\calo_{V^{(d-1)}}[z].$$

We may assume, in addition, that a sequence of monomial transformation has been defined so that the singularities are in the monomial case. Under these assumptions, a prominent role is played by the constant term $a_{p^{e}}$ in the study of invariants in Section  \ref{sec3}, whereas the information of the other intermediate coefficients  is somehow encoded by other $d-1$-dimensional invariants. The protagonist role of the constant term, that appears for example in formula (\ref{slope}), seems to resemble a reduction to the study of equations of the form
$$f_{p^{e}}(z_1)=z^{p^{e}}+a_{p^{e}}^{}\in\calo_{V^{(d-1)}}[z],$$
also known as equations in the purely inseparable case.

These particular features will appear in Section \ref{sec4}, were it is shown how the techniques discussed in the previous sections lead to resolution of two dimensional schemes.

The first proof of resolution in dimension 2, is due to Abhyankar (see \cite{Ab3}). It gives a non-embedded proof of resolution of surfaces in positive characteristic. A synthetic and detailed presentation of this proof appears also in \cite{Cut3}.
Resolution of arithmetical surfaces (non-equicharacteristic surfaces) was proved by J. Lipman making use of techniques of duality theory (\cite{Lip}). 

 Hironaka proves embedded resolution of 2-dimensional schemes over algebraically closed fields in \cite{Hir67}. Moreover, such resolution is attained, as his general proof in characteristic zero, by successive blow-ups along centers included in the Hilbert-Samuel stratum.
More recently, in \cite{CJS}, embedded resolution, in the sense of Hironaka, has been proved for 2-dimensional schemes in the non-equicharacteristic case, namely for arithmetical schemes.

So the outcome of Section \ref{sec4}, namely that of embedded resolution of two dimensional schemes, is to be taken simply as an application of the alternative form of induction discussed here. 
In this last section, the resolution of surfaces is divided in five sub-cases: A), B), C), D1), D2). Each case is defined in accordance to the values attained by the $d-1$-dimensional invariants discussed in Section \ref{sec3}.

In studying the effect of quadratic transformations within case A) we make use of an invariant which appears in Hironaka's proof in \cite{Hir67}, expressed there in terms of Newton polygons, and known as the $\beta$-invariant of the singularity. A similar notion appears also in Abhyankar's work in \cite{Ab3} (see Abhyankar's trick in Lemma 4.12). This invariant is also considered in \cite{Co} and \cite{HaWa}.

The task, for future research, would be to apply these techniques to the open problem
of embedded resolution of three dimensional schemes. 
Abhyankar proves resolution in dimension 3, over algebraically closed fields of characteristic $p>5$ in \cite{Ab3d}. This is a non-embedded proof, that make use of geometric arguments (see also \cite{Cut}).
More recently Cossart and Piltant in \cite{CP1} prove a theorem of resolution in dimension 3, which holds in any characteristic. The outcome of their procedure is very strong, and  it is close to that of embedded resolution, as their resolution only modifies singular points. This last result suggests that embedded resolution of singularities should hold, at least for 3-dimensional schemes over perfect fields.
\end{parrafo}

\begin{parrafo}{\bf Acknowledgement}. We have profited from discussions with A. Bravo and V. Cossart. We thank the referee for the comments and corrections that helped us to improve the  presentation of this paper.
\end{parrafo}

\section{Rees algebras, elimination and monoidal transformations}\label{sec2}

\begin{parrafo}
In problems concerning resolution of singularities it is natural to consider data with two ingredients. The first ingredient is given by a hypersurface, say $X$, embedded in a $d$-dimensional smooth scheme, say $V^{(d)}$, and the second is a positive integer, say $b$. 

A first motivation for this  approach appears already when we fix $b$ as the highest multiplicity of  $X$. In such case one considers the data $(I(X),b)$, where $I(X)\subset\calo_{V^{(d)}}$ is the ideal of definition of $X$. This, in turn, defines a closed set:
$$\{x\in V^{(d)}\ |\ \nu_x(I(X))\geq b\},$$
where $\nu_x(I(X))$ is the order of $I(X)$ at the regular local ring $\calo_{V^{(d)},x}$.
Namely, the set of points where $X$ takes the highest multiplicity $b$ (i.e., the set of $b$-fold points of $X$). 

If $Y$ is a smooth closed center included in the set of $b$-fold points, it  defines a ``transform'' of  the data $(I(X),b)$. In this case the transform is given by the strict transform  of the hypersurface, which makes use of the integer coordinate $b$. In fact, if $V^{(d)}\longleftarrow V_1^{(d)}$ is the blow-up at $Y$, and if $X_1$ is the strict transform of $X$, then
$$I(X)\calo_{V^{(d)}_1}=I(H)^bI(X_1),$$
where $H\subset V_1^{(d)}$ denotes the exceptional hypersurface. The new data $(I(X_1),b)$ will be called the ``transform'' of $(I(X),b)$. If the closed set attached to the new data is empty, $X_1$ has no $b$-fold points and we shall say that the transformation defines a \emph{resolution of $(I(X),b)$}. In this case we have come closer  to the embedded desingularization of the hypersurface.

Hironaka reduces the problem of desingularization to a simultaneous treatment of data of the form, say $(I(X_1),b_1),\dots,(I(X_s),b_s)$. More precisely, to a simultaneous resolution of these previous data by means of monoidal transformations.

In this work we encode these previous data in algebraic terms. This will lead us to a reformulation of resolution in terms of Rees algebras. The Rees algebra attached to the data $(I(X_1),b_1),\dots,(I(X_s),b_s)$ will be the $\calo_{V^{(d)}}$-algebra of the form $\calo_{V^{(d)}}[I(X_1)W^{b_1},\dots, I(X_s)W^{b_s}]$, as we explain below.
\end{parrafo}

\begin{parrafo}
A \emph{Rees algebra} over $V^{(d)}$ is an algebra of the form $\G=\bigoplus_{n\in\N} I_nW^n$, where $I_0=\calo_{V^{(d)}}$ and each $I_n$ is a coherent sheaf of ideals. Here $W$ denotes a dummy variable introduced to keep track of the degree, so $\G\subset\calo_{V^{(d)}}[W]$ is an inclusion of graded algebras. It is always assumed that, locally at any point of $V^{(d)}$, $\G$ is a finitely generated $\calo_{V^{(d)}}$-algebra.  Namely, that the restriction of $\G$ to an affine set $U\subset V^{(d)}$ is of the form
$${\mathcal G}=\calo_{V^{(d)}}(U)[f_{b_1}W^{b_1}, \dots ,f_{b_s}W^{b_s} ](\subset  \calo_{V^{(d)}}[W])$$
for some \emph{local generators}  $\{f_{b_1},\dots,f_{b_s}\}$, each $f_{b_i}\in\calo_{V^{(d)}}$.

We now set the \emph{singular locus of} $\G=\bigoplus I_nW^n$ to be the closed set:
$$\Sing(\G):=\{x\in V^{(d)}\ | \ \nu_x(I_n)\geq n\hbox{ for each }n\in\N\}.$$
It can be checked that for $U$ as before, then
$\Sing(\G)\cap U=\bigcap \{x\in V^{(d)}\ |\ \nu_x(\id{f_{b_i}})\geq b_i\}$.

Fix a monoidal transformation $V^{(d)}\overset{\pi_C}{\longleftarrow}V^{(d)}_1$  with  center $C\subset\Sing(\G)$. For all $n\in\nat$, $I_n\calo_{V_1^{(d)}}$ admits a factorization of the form
$$I_n\calo_{V_1^{(d)}}=I(H_1)^n\cdot I_n^{(1)},$$
where $H_1=\pi_C^{-1}(C)$ denotes the exceptional hypersurface, and $I(H_1)$ the ideal defining $H_1$. 
This defines a Rees algebra over  $ V_{1}^{(d)}$, namely $\G_1=\bigoplus_{n\in\nat}I_n^{(1)}W^n$, called the \emph{transform of $\G$}, denoted by
\begin{equation}\label{trRA}
\xymatrix@R=0pc@C=0pc{
\G & & & & & \G_1\\
V^{(d)}  &  & & & &   V_{1}^{(d)}\ar[lllll]_{\pi_C}.
}\end{equation}

A sequence of transformations will be denoted by:
\begin{equation}\label{seqintro1}
\xymatrix@R=0pc@C=0pc{
\G & & & & & \G_1 &  & & & &  &  & & & &   \G_r\\
V^{(d)}  &  & & & &   V_{1}^{(d)}\ar[lllll]_{\pi_{C}}   & & & & & \dots \ar[lllll]_{\pi_{C_1}} &  & & & &   V_{r}^{(d)}\ar[lllll]_{\pi_{C_{r-1}}}  \\
}
\end{equation}
and herein we always assume that the exceptional locus of the composite morphism $V^{(d)}  \longleftarrow V_r^{(d)}$, say $\{H_1,\dots,H_r\}$, is a union of hypersurfaces with only normal crossings in $V^{(d)}_r$.

A sequence (\ref{seqintro1}) is said to be a \emph{resolution of} $\G$ if, in addition, $\Sing(\G_r)=\emptyset$.


 \end{parrafo}

\begin{parrafo}{\bf Hironaka's main invariants}\label{tauG}. 

We shall circumvent the precise definition of Hironaka's notion of \emph{invariant} (see \cite[2.8]{BeV2}, and \ref{par312} below). But let us mention that his  notion of invariant at a point $x\in\Sing(\G)$ relates to the local codimension of $\Sing(\G)$ at $x$, and also to the codimension of $\Sing(\G_r)$ at points lying over $x$, for sequences as (\ref{seqintro1}). 

There are two main invariants, introduced by Hironaka, which play a crucial role in his Theorem of resolution of singularities (\cite{Hir64}). Both are characteristic free, and we shall formulate them within the context of Rees algebras.

(1) \emph{Hironaka's $d$-dimensional function}. Fix a Rees algebra $\G=\oplus_n I_nW^n$, Hironaka's $d$-dimensional \emph{order} function, say
$$\ord:\Sing(\G)\longrightarrow \mathbb{Q}_{>0}$$
is defined by setting
$$\ord(\G)(x)=\min_n\Big\{\frac{\nu_x(I_n)}{n}\Big\},$$
where $\nu_x(I_n)$ denotes the order of $I_n$ at the regular local ring $\calo_{V^{(d)},x}$. Here the invariant at $x\in\Sing(\G)$, in Hironaka's sense, is the value $\ord(\G)(x)$.

\vspace{0.2cm}

(2) \emph{Hironaka's $\tau$-invariant}. This is a positive integer attached here  to every closed point $x\in\Sing(\G)$. Recall that the tangent space at $x$ is
${\mathbb T}_{V^{(d)},x}=\Spec(gr_M(\calo_{V^{(d)},x}))$, where $gr_M(\calo_{V^{(d)},x})$ is the graded ring of the regular local ring $\calo_{V^{(d)},x}$.
An homogeneous ideal $\In_x(\G)$ in $gr_M(\calo_{V^{(d)},x})$ is defined by $\G=\oplus_n I_nW^n
$ at any $x\in \Sing(\G)\subset V^{(d)}$.  $\In_x(\G)$ is the ideal spanned by the class of $I_n$ in $M^n/M^{n+1}$, for all $n\geq 1$.

The \emph{tangent cone of $\G$}, say $\mathcal C_{\G,x}\subset \mathbb T_{V^{(d)},x},$ is the homogeneous subscheme defined by $\In_x(\G)$ in $gr_M(\calo_{V^{(d)},x})$.

Here we view the tangent space as a vector space.  A subspace $S\subset \mathbb T_{V^{(d)},x}$ acts by translation defining, in this way, an additive group scheme over $k(x)$. More precisely, set $tr_v(u)=u+v$ with $v\in S$ and $u\in{\mathbb T}_{V^{(d)},x}$.
Define $\L_{\G,x}\subset \mathcal C_{\G,x}$ to be the biggest additive subscheme so that $\C_{\G,x}+\L_{\G,x}=\C_{\G,x}$ (the biggest subspace acting on the tangent cone of $\G$). This is called the \emph{subscheme of vertices} of $ \mathcal C_{\G,x}$. 

Finally, define the \emph{$\tau$-invariant at $x$}, say $\tau_{\G,x}$, as the codimension of 
 $\L_{\G,x}$ in ${\mathbb T}_{V^{(d)},x}$.
 
 \end{parrafo}

\begin{parrafo}\label{rp23}{\bf Rees algebras and differential structure}. 

There is a curious \emph{compatibility} of differential operators on smooth schemes and Hironaka's notion of invariants. Let $V^{(d)}$ denotes a smooth scheme over $k$. 

A Rees algebra  ${\mathcal G}=\bigoplus_{n}I_nW^n$ over  $V^{(d)}$ is said to be a {\em differential Rees algebra over $k$} if taking restrictions of $\G$ over  every open affine set, 
$D_{r}(I_n)\subset I_{n-r},$
 for any index $n$ and for any $k$-differential operator $D_r$  of order $r<n$. When a smooth morphism  of $k$-schemes, say $V^{(d)}\overset{\beta}{\longrightarrow}  V^{(d')}$, is fixed, and the previous property holds 
for differential operators which are $\calo_{V^{(d')}}$-linear, or say, 
$\beta$-relative operators, then $\G$ is said to be a \emph{$\beta$-relative differential Rees algebra}, or simply \emph{$\beta$-differential}.

\begin{proposition}$($\cite[Theorems 3.2 and 4.1]{VV1}$)$. 
Every Rees algebra $\G$ over $V^{(d)}$ admits an extension to a new Rees algebra, say $\G\subset Diff(\G)$, so that $Diff(\G)$ is a differential Rees algebra. Moreover, this differential algebra has the following properties:
\begin{enumerate}
\item $Diff(\G)$ is the smallest differential Rees algebra containing $\G$.

\item $\Sing(\G)=\Sing(Diff(\G))$.

\item The equality in $(2)$ is preserved by transformations. In particular, any resolution of $\G$ defines a resolution of $Diff(\G)$, and the converse also holds.
\end{enumerate}
\end{proposition}

The property in (3) says that, for the sake of defining a resolution of $\G$, we may always assume that it is a differential Rees algebra. This is an important reduction because, as we shall indicate, differential Rees algebras have very handy properties. In fact, they turn out to be very useful when defining projections and other structures introduced by these projections, as we discuss below.
\end{parrafo}

\begin{parrafo}{\bf Transversal projections and elimination}\label{trans}.
Once we fix a closed point $x\in V^{(d)}$ it is very simple to construct, for any positive integer $d'\leq d$, a smooth scheme $V^{(d')}$ together with a smooth morphism $\beta:V^{(d)}\longrightarrow V^{(d')}$ (a projection); at least after restriction of $V^{(d)}$ to an \'etale neighborhood of $x$. This claim follows, essentially, from the fact that $(V^{(d)},x)$ is an \'etale neighborhood of the affine space $\A^{(d)}$ at the origin, say $(\mathbb{A}^{(d)},\mathbb{O})$ (see \cite{AlKe}). Plenty of smooth morphisms between affine spaces can be constructed (in fact, plenty of surjective linear maps), say $(\mathbb{A}^{(d)},\mathbb{O})\longrightarrow (\mathbb{A}^{(d')},\mathbb{O})$, for $d'\leq d$. 
In fact, the smoothness of $V^{(d)}$ over the perfect field $k$ ensures that if $\{x_1,\dots,x_d\}$ is a regular system of parameters at $\calo_{V^{(d)},x}$, then $(V^{(d)},x)$ is an \'etale neighborhood of $\mathbb{A}^{(d)}_k=\Spec(k[x_1,\dots,x_d])$ at the origin.

Furthermore, whenever we fix a subspace $S$ of dimension $d-d'$ in  ${\mathbb T}_{V^{(d)},x}$, a smooth scheme
$ V^{(d')}$, together with a smooth morphism $\beta:V^{(d)}\longrightarrow V^{(d')}$, can be constructed so that 
$ker (d(\beta)_x)=S$ (here $d(\beta)_x:{\mathbb T}_{V^{(d)},x} \longrightarrow {\mathbb T}_{V^{(d')},\beta(x)}$ is a surjective linear transformation).

Fix now a differential Rees algebra over $V^{(d)}$ and a closed point $x\in\Sing(\G)$. Set $\tau_{\G,x}=e$ (the codimension of $\mathcal{L}_{\G,x}(\subset\mathcal{C}_{\G,x})$ in the tangent space $\mathbb{T}_{V^{(d)},x}$). 

Fix $d'$ so that $d \geq d'\geq d-e$. We say that a smooth morphism  $\beta:V^{(d)}\longrightarrow V^{(d')}$ is {\em transversal} to $\G$ at $x$ if 
$$ker(d\beta)_x\cap\mathcal{L}_{\G,x}=\mathbb{O}.$$ 
Here, $ker(d\beta)_x$ has dimension $=d-d'\leq e$ and the previous condition says that both spaces are in general position. Moreover, this is an open condition  which  holds at points in an open neighborhood of $x$ (see \cite[Remark 8.5]{BV3}).

\begin{definition}\label{def26}
A smooth morphism  $\beta:V^{(d)}\longrightarrow V^{(d')}$ is said to be {\em transversal} to $\G$, if 
\begin{enumerate}
\item $\tau_{\G,x}\geq d-d'$, and
\item $ker(d\beta)_x\cap\mathcal{L}_{\G,x}=\mathbb{O}$,
\end{enumerate}
at any closed point $x\in\Sing(\G)$.
\end{definition}

Notice that if $\G$ is a differential Rees algebra, then it is, in particular, a $\beta$-differential Rees algebra for \emph{any} transversal morphism $\beta$. The usefulness of differential Rees algebras relies on this particular fact, as we shall see in the following proposition.

\begin{proposition}\label{elim} Assume that $\beta:V^{(d)}\longrightarrow V^{(d')}$ is transversal to $\G$, and that $\G$ is $\beta$-differential.
Then a Rees algebra, say $\R_{\G,\beta}$, is defined 
over $V^{(d')}$,
\begin{equation}\label{lec1}
\xymatrix@C=2.5pc@R=-0.15pc{
\G & \R_{\G,\beta}\\
V^{(d)}\ar[r]^{\beta} & V^{(d')},
}
\end{equation}
$($i.e., $\R_{\G,\beta}\subset\calo_{V^{(d')}}[W])$. $\R_{\G,\beta}$ is called the {\em elimination algebra of $\G$ defined by $\beta$}, and has the following properties:
\begin{enumerate} 
\item The natural lifting of $\R_{\G,\beta}$ to $\calo_{V^{(d)}}$, say $\beta^*(\R_{\G,\beta})\subset \calo_{V^{(d)}}[W]$, is a subalgebra of $\G$. $($\cite[Theorem 4.13]{VV4}$)$.

\item $\beta({\Sing(\G)}) \subset \Sing(\R_{\G,\beta})$, and
$\beta|_{\Sing(\G)}:\Sing(\G)\longrightarrow \Sing(\R_{\G,\beta})$ defines  a set theoretical bijection of $\Sing(\G)$ with its image.  Moreover, corresponding points of $\Sing(\G) (\subset V^{(d)})$ and $\Sing(\R_{\G,\beta})(\subset V^{(d')})$ have the same residue field. $($\cite[1.15 and Theorem 4.11]{VV4}, or \cite[7.1]{BV3}$)$.

\item Given a smooth sub-scheme $Y\subset\Sing(\G)$, then $\beta(Y)(\subset\Sing(\R_{\G,\beta}))$ is isomorphic to $Y$. In particular $Y$ defines a transformation of $\G$ and also of $\R_{\G,\beta}$. $($\cite[Theorem 9.1 (i)]{BV3}$)$.

\item $($\cite[Theorems 10.1 and 9.1]{BV3}$)$. A smooth center $Y\subset\Sing(\G)$ defines a commutative diagram
\begin{equation}\label{lec2}\xymatrix@R=-0.15pc@C=3pc{
\G &  \G_1\\
V^{(d)}\ar[dddd]^{\beta} & V_1^{(d)}\ar[l]_{\pi_Y}\ar[dddd]^{\beta_1}\\
\\
\\
\\
V^{(d')} & V_1^{(d')}\ar[l]_{\pi_{\beta(Y)}}\\
\R_{\G,\beta} & (\R_{\G,\beta})_1
}\end{equation}
where $\G_1$ and $(\R_{\G,\beta})_1$ denote the transforms of 
$\G$ and $\R_{\G,\beta}$,  respectively.  Here $\beta_1$ is uniquely determined, and defined in the restriction of $V_1^{(d)}$ to a neighborhood of $\Sing(\G_1)$. This diagram has the following additional properties:

{\rm(4a)} $V_1^{(d)}\overset{\beta_1}{\longrightarrow} V_1^{(d')}$ is transversal to 
$\G_1$ and $\G_1$ is $\beta_1$-differential. In particular, we get:
$$\xymatrix@C=2.5pc@R=0pc{
\G_1 & \R_{\G_1,\beta_1}\\
V_1^{(d)}\ar[r]^{\beta_1} & V_1^{(d')}
}$$

{\rm (4b)}   $(\R_{\G,\beta})_1$ coincides with $\R_{\G_1,\beta_1}$, the elimination algebra 
 of $\G_1$ defined by $\beta_1$.

\item $($\cite[Theorem 10.1]{BV3}$)$ If a different transversal morphism, say $\beta':V^{(d)}\longrightarrow \widetilde{V}^{(d')}$, defining
$$\xymatrix@C=2.5pc@R=-0.15pc{
\G & \R_{\G,\beta'}\\
V^{(d)}\ar[r]^{\beta'} & \widetilde{V}^{(d')},
}
$$
is considered, then analogous properties to $(1)$, $(2)$, $(3)$ and $(4)$ hold. Moreover, the order of both elimination algebras coincide at any point $x\in\Sing(\G)$, i.e.,
$$\ord(\R_{\G,\beta})(\beta(x))=\ord(\R_{\G,\beta'})(\beta'(x)).$$
This says a function
$$\Ord^{(d')}(\G):\Sing(\G)\longrightarrow \mathbb{Q}_{>0}$$
can be defined by setting 
$$\Ord^{(d')}(\G)(x)=\ord(\R_{\G,\beta})(\beta(x)).$$
\end{enumerate}
\end{proposition}

 

\end{parrafo}

\begin{parrafo}\label{cuad}
The previous proposition states that given 
 $\beta:V^{(d)}\longrightarrow V^{(d')}$ transversal to $\G$, and assuming that 
 $\G$ is a $\beta$-differential Rees algebra (e.g., $\G$ is an absolute differential Rees algebra), then
a sequence of  transformations defined with the conditions in (\ref{seqintro1}), say
\begin{equation}\label{cmut}\xymatrix@R=0pc{
\G & \G_1  & & \G_r\\
V^{(d)} & V_1^{(d)}\ar[l]_{\pi_{Y}} & \dots\ar[l]_{\pi_{Y_1}} &V_r^{(d)}\ar[l]_{\pi_{Y_{r-1}}}
}\end{equation}
gives rise to a diagram, say
\begin{equation}\label{conmut}
\xymatrix@R=0pc@C=3.5pc{
\G & \G_1  & & \G_r\\
V^{(d)}\ar[dddd]_ {\beta} & V_1^{(d)}\ar[l]_{\pi_{Y}}\ar[dddd]_ {\beta_1} & \dots\ar[l]_{\pi_{Y_1}} &V_r^{(d)}\ar[dddd]_ {\beta_r}\ar[l]_{\pi_{Y_{r-1}}}\\
\\
\\
\\
V^{(d')} & V_1^{(d')}\ar[l]_{\pi_{\beta(Y)}} & \dots\ar[l]_{\pi_{\beta_1(Y_1)}} &V_r^{(d')}\ar[l]_{\pi_{\beta_{r-1}(Y_{r-1})}}\\
\R_{\G,\beta} & (\R_{\G,\beta})_1 & & (\R_{\G,\beta})_r
}\end{equation}
where:
\begin{enumerate}
\item For any index $i$, there is an inclusion $(\R_{\G,\beta})_i\subset \G_i$.

\item  $\beta_i(\Sing(\G_i))\subset \Sing((\R_{\G,\beta})_i)$, and $\beta_i|_{\Sing(\G_i)}:\Sing(\G_i)\longrightarrow \beta_i(\Sing(\G_i))$ is an identification.

\item Every $\beta_i$ is transversal to $\G_i$ and $\G_i$ is $\beta_i$-differential.

\item $\xymatrix@C=2.5pc{V_i^{(d')} & V^{(d')}_{i+1}\ar[l]_{\pi_{\beta(Y_i)}}}$ denotes the transformation with center $\beta_i(Y_i)$ (isomorphic to $Y_i$).

\item $(\R_{\G,\beta})_i=\R_{\G_i,\beta_i}$ where the later denotes the elimination algebra of $\G_i$ with respect to $\beta_i$.
\end{enumerate}
In the characteristic zero case, the inclusions in (2) are equalities, but in positive characteristic, in general, only the inclusions are ensured.
\end{parrafo}

\begin{remark} Given a Rees algebra $\G$ over $V^{(d)}$, the aim is to construct a sequence as (\ref{cmut}) which defines a resolution.  For this purpose, we can take  $\G$  to be a differential Rees algebras. This additional condition ensures that, whenever we construct a transversal morphism $\beta:V^{(d)}\longrightarrow V^{(d')}$, $\G$ will be a $\beta$-differential Rees algebra. In fact an absolute differential Rees algebra is always relative differential.

Recall that, passing from $\G$ to $Diff(\G)$ does not affect the singular locus, i.e., $\Sing(\G)=\Sing(Diff(\G))$, and moreover the $\tau$-invariant does not change, i.e., $\tau_{\G,x}=\tau_{Diff(\G),x}$ at any closed point $x\in\Sing(\G)$.

The $\tau$-invariant has very subtle implications in resolution problems. For one thing, $\tau_{\G,x}$ is an upper bound of the local codimension of the closed set $\Sing(\G)\subset V^{(d)}$ at the point $x$. If equalities holds, i.e., $\tau_{\G,x}=codim_x(\Sing(\G)$, then $\Sing(\G)$ is smooth in a neighborhood of $x$, and the resolution can be achieved by blowing-up at $\Sing(\G)$. In particular, if $\tau_{\G,x}=d$, then $\Sing(\G)=\{x\}$ (an isolated point) and the quadratic transformation at $x$ defines a resolution of $\G$. Therefore, the strategy is to define resolution of Rees algebras on $V^{(d)}$ by decreasing induction on the highest value of $\tau$.

The following theorem is stated under the inductive assumption of existence of resolution of Rees algebras $\G'$ which fulfills the condition $\tau_{\G',x}\geq d-d'+1$ at any closed point $x\in\Sing(\G')$. 
\end{remark}

\begin{theorem}$($\cite[10.4]{BV3}$)$.\label{thm:BV}
Let $\G$ be a differential Rees algebra. Assume that $\tau_{\G,x}\geq d-d'$ at any closed point $x\in\Sing(\G)$. There is a sequence of  transformations $(\ref{cmut})$,
so that for any local transversal projection 
 $\beta:V^{(d)}\longrightarrow V^{(d')}$, the induced sequence in the lower row of $(\ref{conmut})$ is either a resolution, or is such  that $(\R_{\G,\beta})_r$ is a monomial algebra supported on the exceptional locus. Furthermore, in the latter case the monomial algebra $\beta^{*}_r((\R_{\G,\beta})_r)$ is independent of $\beta$.
\end{theorem}

\begin{parrafo} In the case of characteristic zero, it is easy to extend the sequence in Theorem \ref{thm:BV} to a resolution. In fact, in such a case, this extension is obtain by blowing-up centers prescribed by the monomial algebra $(\R_{\G,\beta})_r$. However, this latter simple construction, that grows from the fact that 
 $\beta_i(\Sing(\G_i))=\Sing((\R_{\G,\beta})_i)$ in (\ref{conmut}),(2), does not apply in positive characteristic. 

In the forthcoming Section \ref{sec3}, we study numerical conditions under which a sequence (\ref{cmut}) in the conditions of Theorem \ref{thm:BV} does extend to a resolution (Theorem \ref{thm33}). In Section \ref{sec4} we show that these numerical conditions can be easily achieved in low dimension. This fact  enables us to prove resolution of singularities of $2$-dimensional schemes over perfect fields.

\end{parrafo}

\section{H-functions, tight monomial algebras and presentations}\label{sec3}

\begin{parrafo}Fix a Rees algebra $\G$ over a $d$-dimensional smooth scheme $V^{(d)}$. At any closed point $x\in\Sing(\G)$, the invariant $\tau_{\G,x}$ was defined as the codimension of the subspace $\mathcal{L}_{\G,x}\subset\mathbb{T}_{V^{(d)},x}$. The requirement that $\tau_{\G,x}\geq d-d'$  at every closed point $x\in\Sing(\G)$, ensures that the conditions in Definition \ref{def26} holds for a generic smooth morphism $\beta:V^{(d)} \longrightarrow V^{(d')}$. 

Under the previous conditions, a so called {\em $d'$-dimensional H-function}, say
$$\Hord^{(d')}(\G):\Sing(\G)\longrightarrow\mathbb{Q}_{>0},$$
was defined in \cite{BeV2}. The definition of this function involves a transversal morphism $\beta$, and the induced elimination algebra $\R_{\G,\beta}$, among other things. Further details will be outlined below (see \ref{explicit}). When specializes to the case of characteristic zero, this function coincides with the \emph{inductive function} introduced by Hironaka, which turned out to be extremely useful in resolution theorems.


Set now a sequence of transformations of $\G$, say
\begin{equation}\label{cmut2}
\xymatrix@R=0pc@C=3pc{
\G & \G_1 & & \G_r\\
V^{(d)} & V_{1}^{(d)}\ar[l]_{\pi_Y} &\dots\ar[l]_{\ \ \ \pi_{Y_1}} & V_r^{(d)}\ar[l]_{ \ \pi_{Y_{r-1}}}.
}\end{equation}
At any step of the sequence, $i=1,\dots,r$, functions
\begin{equation}\label{Hi}
\Hord_i^{(d')}(\G_i):\Sing(\G_i)\longrightarrow \mathbb{Q}_{>0}
\end{equation}
are defined.

Moreover, once we fix a transversal projection $\beta:V^{(d)}\longrightarrow V^{(d-1)}$ as above, the previous sequence gives rise to a diagram of the form
\begin{equation}\label{conmut2}
\xymatrix@R=0pc@C=3.5pc{
\G & \G_1  & & \G_r\\
V^{(d)}\ar[dddd]_ {\beta} & V_1^{(d)}\ar[l]_{\pi_{Y}}\ar[dddd]_ {\beta_1} & \dots\ar[l]_{\pi_{Y_1}} &V_r^{(d)}\ar[dddd]_ {\beta_r}\ar[l]_{\pi_{Y_{r-1}}}\\
\\
\\
\\
V^{(d')} & V_1^{(d')}\ar[l]_{\pi_{\beta(Y)}} & \dots\ar[l]_{\pi_{\beta_1(Y_1)}} &V_r^{(d')}\ar[l]_{\pi_{\beta_{r-1}(Y_{r-1})}}\\
\R_{\G,\beta} & (\R_{\G,\beta})_1 & & (\R_{\G,\beta})_r
}
\end{equation}
with the five properties of (\ref{conmut}). In this setting, given $x_i\in\Sing(\G_i)$ the value $\Hord_i(d')(\G_i)(x_i)$ can be computed from $\beta$, and from the liftings $\beta_i$ (in a explicit manner). Moreover, the explicit calculation of the function will lead to the inequalities:
$$
\Hord_i^{(d')}(\G_i)(x_i)\leq \ord((\R_{\G,\beta})_i)(\beta_i(x_i)).
$$

A particular feature of the H-functions is their unpredictable behavior under blow-ups. Another aspect is that they are not upper semi-continuous.
The interest of the previous inequalities is that the functions in the right hand side, say $\ord((\R_{\G,\beta})_i)$, are upper semi-continuous (i.e., H-functions are upper bounded by upper semi-continuous functions).

Before we formulate some further applications of the H-functions in Theorem \ref{thm33}, let us indicate how they  lead to the definition of a monomial algebra attached to an arbitrary sequence (\ref{cmut2}).

\begin{definition}\label{tight}
Fix a differential Rees algebra $\G$, a transversal projection $\beta:V^{(d)}\longrightarrow V^{(d')}$, and a sequence of transformations (\ref{cmut2}). The \emph{tight monomial algebra} attached to the sequence (\ref{cmut2}) is a monomial algebra supported on the exceptional locus, say
\begin{equation}\label{eqmon}
\m_r W^s=\calo_{V_r^{(d')}}[I(H_1)^{h_1}\dots I(H_r)^{h_r}W^s],
\end{equation}
where the exponents $h_i$ are such that
$$\frac{h_i}{s}=\Hord^{(d-1)}_{i-1}(\G_{i-1})(\xi_{Y_{i-1}})-1.$$
Here $\xi_{Y_{i-1}}$ denotes the generic point of $Y_{i-1}$, the center of the blow-up.
\end{definition}

The real strength of the tight monomial algebra appears in the statement of the following theorem. 

\begin{theorem}\label{thm33} $($\cite[Theorem 6.6]{BeV1}$)$.
Fix a differential Rees algebra $\G$ and a transversal projection $\beta:V^{(d)}\longrightarrow V^{(d')}$. Then the following inequalities hold for any sequence as $(\ref{cmut2})$:
$$\ord(\m_rW^s)(\beta_r(x))\leq \Hord_r^{(d')}(\G_r)(x)\leq \ord((\R_{\G,\beta})_r)(\beta_r(x))$$
and all $x\in\Sing(\G_r)$. That is, the upper semi-continous functions $\ord(\m_rW^s )$ and $\ord((\R_{\G,\beta})_r)$ are lower and upper bounds, respectively, of the function $\Hord^{(d')}_r(\G_r)$.

Moreover, if equality holds in the left-hand side of the previous inequalities,  then the sequence $(\ref{cmut2})$ can be extended to a resolution of $\G$. 
\end{theorem}

The second half of the previous theorem indicates that the resolution of the Rees algebra $\G$ can be achieved if some suitable numerical conditions hold.

In Section \ref{sec4} it is shown that such numerical conditions can be attained if $\G$ is a Rees algebra attached to the resolution of a $2$-dimensional scheme, leading to a resolution of $2$-dimensional schemes \emph{a la Hironaka} (by successively blowing-up along smooth centers included in the highest Hilbert-Samuel stratum).

\end{parrafo}

\begin{parrafo}\label{explicit}{\bf Explicit computation of the H-functions}. 

In our previous discussion we have treated some properties of H-functions. Here we indicate how to compute explicitly the values $\Hord^{(d')}(\G)(x)$ at points $x\in\Sing(\G)$. Further details about this computation can be found in \cite{BeV1} and \cite{BeV2}.

We first address the case $d'=d-1$. So assume $\tau_{\G,x}\geq 1$ for any $x\in\Sing(\G)$. To pave the way to the definition of the H-function in this context,  we define, locally at a point $x\in\Sing(\G)$ a concept of \emph{$p$-presentation of $\G$}, say 
\begin{equation}\label{ppp}
p\P(\beta,  z,  f_{p^e}(z)).
\end{equation}
These data consist of:
\begin{enumerate}
\item a transversal projection $\beta:V^{(d)}\longrightarrow V^{(d-1)}$ (see Definition \ref{def26}),
\item a global function on $V^{(d)}$, say $z$,  so that $\{z=0\}$ defines a section of $\beta$ (i.e., $\{dz\}$ is a basis of the module of $\beta$-differentials, say $\Omega_\beta^1$), and
\item a monic polynomial of order $p^e$, say 
$$f_{p^e}(z)=z^{p^e}+a_1z^{p^e-1}+\dots+a_{p^e}\in \calo_{V^{(d-1)}}[z],$$
where each $a_i$ is a global function on $V^{(d-1)}$, and $p$ is the characteristic of the underlying perfect field $k$.
\end{enumerate}
In \cite[Proposition 2.11]{BeV1} it is proved that such data can always be locally defined, when $\G$ is $\beta$-differential. This last requirement imposed no serious conditions as we may assume that $\G$ is a differential algebra. Moreover, it is shown that  $\G$ and the graded algebra
\begin{equation} \label{eqdpl} 
\calo_{V^{(d)}}[f_{p^e}(z)W^{p^e},\Delta_z^{j}(f_{p^e}(z))W^{p^e-j}]_{1\leq j\leq p^e-1}\odot\beta^*(\R_{\G,\beta})
\end{equation}
have the same integral closure.

The previous $\Delta_z^j$ are  $\beta$-differential operators defined in terms of the Taylor morphism in the following manner:

Consider the morphism of $S$-algebras $Tay: S[Z]\longrightarrow S[Z, T]$,
 defined by setting $Tay(Z)=Z+T$ (Taylor expansion). Here
\begin{equation}\label{lapiz}
Tay(f(Z))= f(Z+T)= \sum \Delta^{r}(f(Z))T^r,
\end{equation}
and these operators $ \Delta^{r}\!: S[Z] \longrightarrow S[Z]$ are defined by this morphism. It is well known that $\{\Delta^{0}, \Delta^{1}, \dots,\Delta^{r}\}$ is a basis of the free module of $S$-differential operators of order $r$. The same applies here for $\calo_{V^{(d-1)}} [z] $: the set $\{\Delta^{0}_z, \Delta^{1}_z, \dots,\Delta^{r}_z\}$ consists of  differential operators of order $r$ over $V^{(d-1)}\times \mathbb A^1$. Moreover, as the smoothness of $\beta$ ensures that $V^{(d)}$ is \'etale over $V^{(d-1)}\times \mathbb A^1$, the previous set also generates $Diff^r_{\beta}$, the $\beta$-linear differential operators of order $r$ (for $\beta:V^{(d)}\longrightarrow V^{(d-1)}$).


\begin{remark}\label{rmk34}
Suppose given a $p$-presentation of $\G$ as in (\ref{ppp}), together with a diagram (\ref{conmut2}) (with $d'=d-1$). Then, it is proved in \cite[6.2]{BeV1} that  there is natural lifting of $p\P$ to $p$-presentations of $\G_i$, say $p\P_i(\beta_i,z_i,f_{p^e}^{(i)}(z_i))$, for $i=1,\dots,r$. Here $\beta_i:V^{(d)}_i\longrightarrow V^{(d-1)}_i$ are the smooth morphisms in (\ref{conmut2}), and $z_i=0$ define a $\beta_i$-transversal section.
\end{remark}

\end{parrafo}

\begin{parrafo}\label{par35} We now give a first step towards the definition of the H-functions in (\ref{Hi}). Fix a $p$-presentation of $\G$, say $\P(\beta,  z,  f_{p^e}(z))$, with  $f_{p^e}(z)=z^{p^e}+a_1z^{p^e-1}+\dots+a_{p^e}$. Define the \emph{slope of $p\P$ at $y\in V^{(d-1)}$} as
\begin{equation}\label{slope}
Sl(p\P)(y)=\min\Big\{\frac{\nu_y(a_{p^e})}{p^e},\ord(\R_{\G,\beta})(y)\Big\}\in\mathbb{Q}_{\geq 0}.
\end{equation}
We now present a definition that will lead us to the precise value of the H-function in Theorem \ref{thm38}: A $p$-presentation is \emph{well-adapted to $\G$ at $y\in V^{(d-1)}$} when
\begin{enumerate}
\item[(i)] $Sl(p\P)(y)=\ord(\R_{\G,\beta})(y)$, or
\item[(ii)] $Sl(p\P)(y)=\frac{\nu_y(a_{p^e})}{p^e}$ and $In_y(a_{p^e})\in gr(\calo_{V^{(d-1)},y})$ is not a $p^e$-th power.
\end{enumerate}

\begin{remark}\label{rmk36}
Fix a $p$-presentation $p\P(\beta,z,f_{p^e})$ locally at a point $x\in\Sing(\G)$. One can easily modify $z$ and $f_{p^e}(z)$ so as to obtain a new $p$-presentation, say $p\P'(\beta,z',f_{p^e}'(z'))$, which is well-adapted at $\beta(x)$ (see \cite[Section 5]{BeV1}).
\end{remark}

\begin{theorem}\label{thm38}
Fix a Rees algebra $\G$. The H-function, say
$$\Hord^{(d-1)}(\G):\Sing(\G)\longrightarrow \mathbb{Q}_{>0},$$
is defined by 
\begin{equation}\label{Hord}
\Hord^{(d-1)}(\G)(x)=Sl(p\P)(\beta(x))
\end{equation}
where $p\P=p\P(\beta,z,f_{p^e}(z))$ is any $p$-presentation of $\G$, which is well-adapted at $\beta(x)$. This rational value is independent of any choice as long as $p\P$ is well-adapted at $\beta(x)$.
\end{theorem}

\begin{remark}
The value $\Hord^{(d-1)}(\G)(x)$ in (\ref{Hord}) is given by the slope of a well adapted $p$-presentation $p\P=p\P(\beta,z,f_{p^e}(z))$. Note that the value of the slope in (\ref{slope}) relies only in the constant coefficient of $f_{p^e}(z)$, say $a_{p^e}$, and the order of the elimination algebra. This indicates that the contribution of the intermediate coefficients of $f_{p^e}(z)$, say $a_i$ ($1\leq i\leq p^e-1$), is somehow encoded by the constant coefficient $a_{p^e}$ and the elimination algebra. 

This resembles an  expected (but never proved) behavior of resolution problems in positive characteristic: the reduction to the so called purely inseparable polynomials, namely those of the form $f_{p^e}(z)=z^{p^e}+a_{p^e}$.
\end{remark}

In the previous discussion it is assumed that if $p\P=p\P(\beta,z,f_{p^e}(z))$ is a presentation, then $\G$ is a relative $\beta$-differential algebra. This last condition is automatically guaranteed as we assume that $\G$ is an absolute differential algebra, or a transform of an absolute differential algebra.

\begin{remark}
Fix a differential Rees algebra $\G$, a sequence of transformations (\ref{cmut}) and a $p$-presentation $p\P$ of $\G$.  Remark \ref{rmk34} enables us to define $p$-presentations of $\G_i$, say $p\P_i=p\P_i(\beta_i,z_i,f_{p^e}^{(i)}(z_i))$.  Remark \ref{rmk36} indicates how to make use of $p\P_i$ to define a function 
$$\Hord^{(d-1)}_i(\G_i):\Sing(\G_i)\longrightarrow \mathbb{Q}_{>0},$$
for each index $i=0,\dots,r$.
\end{remark}
\end{parrafo}

\begin{parrafo}\label{par310}

Now we address the definition of the H-functions in (\ref{Hi}) for  the general case of arbitrary $d'$ ($1\leq d'\leq d$). Set $\ell$ so that $\tau_{\G,x}\geq d-d'=\ell$ at any closed point $x\in\Sing(\G)$. Let $\beta:V^{(d)}\longrightarrow V^{(d')}$ be a transversal projection.
The problem, in this general case, is to define a notion of $p$-presentations, with similar properties as that in (\ref{ppp}) for de case $d'=d-1$ (i.e., $\ell=1$). The extension to this general case is not straightforward, and has been treated in \cite[Theorem 5.12]{BeV2}. There it is proved that, if $\tau_{\G,x}\geq \ell$ at each closed point $x\in\Sing(\G)$, and $\G$ is a differential algebra or a transformation of a differential algebra, then there is a \emph{$p$-presentation}, say
$$p\P=p\P(\beta,z_1,\dots,z_\ell,f_{p^{e_1}}(z_1),\dots,f_{p^{e_\ell}}(z_\ell)),$$
for which $\G$ has the same integral closure as
$$
\calo_{V^{(d)}}[f_{p^{e_i}}(z_i)W^{p^{e_i}},\Delta_{z_i}^{j_i}(f_{p^{e_i}}(z_i))W^{p^{e_i}-j_i}]_{1\leq j_i\leq p^{e_i}-1,\ 1\leq i\leq \ell}\odot\beta^*(\R_{\G,\beta}).
$$
Moreover the $p$-presentation is defined so that
\begin{itemize}
\item $\{z_1=\dots=z_\ell=0\}$ is a section of the smooth morphism  $\beta$ (i.e., $\{dz_1,\dots,dz_\ell\}$ is a basis of $\Omega_\beta^1$), 
\item  $f_{n_i}(z_i)W^{n_i}\in\G_r$ for each index $i=1,\dots,\ell$, and the polynomials are of the form:
\end{itemize}
\begin{equation}\label{pol2}
\begin{array}{l}
f_{p^{e_1}}(z_1)=z_1^{p^{e_1}}+a_1^{(1)}z_1^{p^{e_1}-1}+\dots+a_{p^{e_1}}^{(1)}\in\calo_{V^{(d-\ell)}}[z_1],\\
\ \ \ \vdots\\
f_{p^{e_\ell}}(z_\ell)=z_\ell^{p^{e_\ell}}+a_1^{(\ell)}z_\ell^{p^{e_\ell}-1}+\dots+a_{p^{e_\ell}}^{(\ell)}\in\calo_{V^{(d-\ell)}}][z_\ell],
\end{array}
\end{equation}
Let us stress here that all coefficients are in dimension $d-\ell$, i.e., $a_{j_i}^{(i)}\in\calo_{V^{(d-\ell)}}$. This fact enables us to extend the previous results in \ref{explicit} (the definitions of slope and  of a presentation well-adapted at a point). This will enable us to compute the values of the H-functions in the general setting. Namely to the extension of Theorem \ref{thm38}.

Firstly fix a $p$-presentation $p\P=p\P(\beta,z_1,\dots,z_\ell,f_{p^{e_1}}(z_1),\dots,f_{p^{e_\ell}}(z_\ell)),$ with $f_{p^{e_i}}(z_i)\in\calo_{V^{(d')}}[z_i]$, as in (\ref{slope}). Let the \emph{slope of $p\P$ at $y\in V^{(d')}$} be
$$Sl(p\P)(y)=\min_{ 1\leq i\leq \ell}\Big\{\frac{\nu_{y}(a_{p^{e_i}}^{(i)})}{p^{e_i}},\ord(\R_{\G,\beta})(y)\Big\}\in\mathbb{Q}_{\geq 0}.$$
Now, we define a $p\P$ presentation to be \emph{well-adapted to $\G$ at $y\in V^{(d')}$} when
\begin{enumerate}
\item[(i)] $Sl(p\P)(y)=\ord(\R_{\G,\beta})(y)$, or
\item[(ii)] $Sl(p\P)(y)=\frac{\nu_y(a_{p^{e_i}})}{p^{e_i}}<\ord(\R_{\G,\beta})(y)$ and $In_y(a_{p^{e_i}})$ is not a $p^{e_i}$-th power,  for some $i\in\{1,\dots,\ell\}$.
\end{enumerate}

Also Remark \ref{rmk36} extends to this setting: Once we fix $p\P$ locally at a point $x\in\Sing(\G)$, then it can be modified into a new one which is well-adapted at $\beta(x)$ (see \ref{par35}).

\begin{theorem}\label{Hordd} Fix a Rees algebra $\G$ and assume that $\tau_{\G,x}\geq d-d'=\ell$ for any $x\in\Sing(\G)$. The H-function
$$\Hord^{(d')}:\Sing(\G)\longrightarrow \mathbb{Q}_{> 0},$$
is defined by
$$\Hord^{(d')}(\G)(x)=Sl(p\P)(\beta(x))=\min_{ 1\leq i\leq \ell}\Big\{\frac{\nu_{\beta(x)}(a_{p^{e_i}}^{(i)})}{p^{e_i}},\ord(\R_{\G,\beta})(\beta(x))\Big\},$$
where now $p\P=p\P(\beta:V^{(d)}\longrightarrow V^{(d')},z_1,\dots,z_\ell,f_{p^{e_1}},\dots,f_{p^{e_\ell}})$ is a $p$-presentation which is assumed to be well-adapted at $\beta(x)$. The value $\Hord^{(d')}(\G)(x)$ is independent of any choice (i.e., independent of the $p$-presentation as long it is well-adapted at $\beta(x)$).
\end{theorem}

\begin{proof}
See \cite[Theorem 5.12 and Definition 5.13]{BeV2}.
\end{proof}

\begin{parrafo}\label{par312}
The previous Theorem indicates how to compute the values of the H-function along $\Sing(\G)$. In fact, given a $p$-presentation of $\G$, it is not difficult to modify it so that at a given point $x\in\Sing(\G)$, it is well-adapted at $\beta(x)$. The Theorem says, of course, that the value $\Hord^{(d')}(\G)(x)$ is independent of the choice of the $p$-presentation.

In \cite{BeV2} further properties of these functions are studied, which are related to Hironaka's notion of invariant. This notion has a very precise meaning in the context of resolution of singularities. Roughly speaking, once $V^{(d)}$ is fixed, Hironaka defines an equivalent relation in the class of Rees algebras over $V^{(d)}$. This relation is defined so that if $\G_1$ and $\G_2$ are equivalent, then $\Sing(\G_1)=\Sing(\G_2)$, and given $x\in\Sing(\G_1)$, also $\tau_{\G_1,x}=\tau_{\G_2,x}$. In \cite{BeV2} it is proved that $\Hord^{(d')}(\G_1)(x)=\Hord^{(d')}(\G_2)(x)$ for any point $x\in\Sing(\G_1)$. Namely, that the value $\Hord^{(d')}(\G_1)(x)$ is an \emph{invariant}. 

This last fact will ultimately guarantee that the functions $\Hord^{(d')}(\G)$ will be useful for the problem of resolution of singularities. In fact, resolution of singularities reduces to resolution of Rees algebras, but only if the latter is defined in a way that two equivalent Rees algebras undergo the same resolutions. 
\end{parrafo}

\end{parrafo}


\section{Embedded resolution of $2$-dimensional schemes}\label{sec4}

Here we address the proof of embedded resolution of $2$-dimensional schemes. We show here how our invariants lead to the resolution of a hypersurface embedded in a $3$-dimensional smooth scheme. The extension of the resolution of the hypersurface case, treated here, to that of arbitrary $2$-dimensional schemes is not straightforward. It is a particular feature of the invariants introduced in the previous sections, studied in \cite{BeV2}, which relies essentially on \ref{par310} and Theorem \ref{Hordd}.

\begin{parrafo}{\bfseries Stratification of the exceptional locus.}\label{par:mred}

We shall take as starting point a diagram (\ref{conmut}), in the setting of Theorem \ref{thm:BV}, for $d=3$. In which case $V_{r}^{(d-1)}=V^{(2)}_r$ is a 2-dimensional smooth scheme. Recall that Theorem \ref{thm:BV} enables us to assume that the elimination algebra over $V^{(2)}_r$ is monomial. Namely, that 
$$(\R_{\G,\beta})_r=I(H_1)^{\alpha_1}\dots I(H_r)^{\alpha_r}W^s,$$
for some integers $\alpha_i\geq 0$.
Assume, in addition, that the tight monomial algebra,  defined in (\ref{eqmon}), is of the form 
\begin{equation}\label{arbol}
\m_r W^s=I(H_1)^{h_1}\dots I(H_r)^{h_r}W^s\ \ \ \quad\hbox{ with }\ 0\leq h_i<s.
\end{equation}
We shall indicate below that the condition $0\leq h_i<s$ can be achieved by blowing-up permissible centers of codimension $2$ in the smooth $3$-dimensional scheme $V^{(3)}_r$. In this case, the tight monomial algebra $\m_r W^s$ subject to this condition is said to be \emph{reduced}. 

This last assumption guarantees that $\Sing(\m_rW^s)\subset V_r^{(2)}$ has no components of codimension $1$, and we claim that neither does $\beta_r(\Sing(\G_r))$, which is, therefore, a finite set of closed points.

Throughout  this section, we will always assume that the tight monomial algebra is reduced; that is, every  quadratic transformation will be followed by a finite sequence of blow-ups at centers of codimension $2$, so as to guarantee that the new tight monomial algebra is reduced. 
 
Locally at any $x\in\Sing(\G_r)$ one can fix a $p$-presentation $p\P(\beta_r,z,f_{p^e}(z))$ well-adapted to $\G_r$ at $\beta_r(x)$,
here, $f_{p^e}(z)=z^{p^e}+a_1z^{p^e-1}+\dots+a_{p^e}\in\calo_{V_r^{(2)}}[z]$. 

The $p$-presentation $p\P$ can be chosen so as to have the following two properties (see \cite[Proposition 5.9 and Definition 3.10]{BeV1}):
\begin{itemize}
\item $(\R_{\G,\beta})_r\subset \m_rW^s$, and
\item  $a_jW^j \in\m_rW^s$ for any $j=1,\dots,p^e$.
\end{itemize}
This already proves the claim: if $h_i \geq s$ for some $H_i$ containing $x$, then $\langle z, I(H_i) \rangle$ defines a smooth 
component of $\Sing(\G_r)$ of codimension 2 in $V^{(3)}_r$. This monoidal transformation induced the identity map over $V^{(2)}_r$.
However $h_i$ drops to $h_i-s$  after blowing up at such component.

So we will assume here that $\Sing(\G_r)$ is a finite set of closed points, and we fix the notation $\x=\beta_r(x)$ along this section.

Recall that 
$$\Hord^{(d-1)}(\G_r)(x)=\min\big\{\frac{\nu_\x(a_{p^e})}{p^e},\ord((\R_{\G,\beta})_r)(\x)\big\}.$$
The previous two properties ensures that 
\begin{equation}\label{eqm}
\ord(\m_rW^s)(\x)\leq\Hord^{(d-1)}(\G_r)(x).
\end{equation}

As $(\R_{\G,\beta})_r$ is, by assumption, monomial, Proposition \ref{elim} (2) ensures that the singular locus of $\G_r$ is entirely included in a union of the exceptional hypersurfaces. We  begin by fixing a suitable stratification of the union of the exceptional hypersurfaces in $V_r^{(2)}$, and then we shall construct a stratification on the successive quadratic transformations defined over $V^{(2)}_r$. Here the quadratic transformations are defined canonically, with centers on the finitely many points of the singular locus. This lead to the following definition:

\begin{definition}
An isolated closed point $x\in\Sing(\G_r)$ is said to be \emph{good} (or green) if 
$$\Hord^{(d-1)}(\G_r)(x)=\ord((\R_{\G,\beta})_r)(\x).$$ 
The point $x$ is said to be \emph{bad} (or purple) if 
$$\Hord^{(d-1)}(\G_r)(x)=\frac{\nu_\x(a_{p^e})}{p^e}<\ord((\R_{\G,\beta})_r)(\x).$$
\end{definition}

The previous definition for a point $x\in\Sing(\G_r)$ leads to a coloring of singular points of $\beta_r(\Sing(\G_r))$ which is a finite set of points in $V^{(2)}_r$. We will say that $\x=\beta_r(x)$ is good (bad) if $x$ is good (bad).
In the same manner, exceptional hypersurfaces will be distinguished in terms of the exponents that arise in $\m_rW^s$ and $(\R_{\G,\beta})_r$:
 
\begin{definition}\label{defcolhyper}
We will say that the exceptional hypersurface $H_i$ is \emph{good} (or green) if $h_i=\alpha_i$ and $H_i$ is \emph{bad} (or purple) if $h_i<\alpha_i$.

\end{definition}

New quadratic transformations over $V_r^{(2)}$ will be defined, this will introduce new exceptional hypersurfaces, each of which will be either good or bad.

\begin{remark}
Note from Definition \ref{tight} that the closed point $\x\in\beta_r(\Sing(\G_r))$ is good (bad) if and only if the exceptional hypersurface introduced by the quadratic transformation at $\x$ is good (bad).
\end{remark}

Our stratification will be defined only along the union of bad hypersurfaces. Singularities lying entirely on the good locus are, to some extend, friendly singularities.
In fact, if a closed point $x$ is only included in good hypersurfaces, then $\ord((\R_{\G,\beta})_r)(\x)=\ord(\m_r W^s)(\x)$ and we claim that the transversal parameter $z$ of the $p$-presentation $p\P$ can be chosen so as to be a hypersurface of maximal contact.  
To this end, in \cite[Theorem 8.5]{BeV1}, using (\ref{eqm}) and the two properties of $\m_rW^s$ mentioned above, it is proved that 
$$\ord((\R_{\G,\beta})_r)(\x)=\ord(\m_r W^s)(\x)\ \Longrightarrow \ (\R_{\G,\beta})_r=\m_r W^s$$ 
in a neighborhood of $\x$.
In such case, as  $a_jW^j\in\m_rW^s$, we conclude that $a_jW^j\in(\R_{\G,\beta})_r\subset \G_r$ for $j=1,\dots,p^e$. Check now that $zW$ fulfills the integral condition 
$\lambda^{p^e}+(a_1W^1)\lambda^{p^e-1}+\dots+(a_{p^e}-f_{p^e}(z))W^{p^e}=0$, which guarantees $zW\in\G_r$. This ensures that $z=0$ defines a hypersurface of maximal contact.
\end{parrafo}

\begin{parrafo}{\bf Stratification at level $r$}.

\vspace{0.1cm}

The image of $\Sing(\G_r)$ is a finite set of points. We will first define a stratification only locally around points of this image.  In \ref{parstr} and \ref{parinf} we shall indicate how to define the stratification after blowing-up these points, and moreover after applying a sequence of quadratic transformations.

According to the previous discussion, we only draw attention to those points of $\Sing(\G_r)$ with images in the union of bad lines of $V^{(2)}_r$. As these points are isolated, we may assume, after restriction, that the point is unique. The following two situations can arise:
\begin{itemize}
\item There exists a unique hypersurface $H_1$ so that $\x\in H_1$ (in which case $H_1$ must be bad).
\item $\x$ is an intersection point of two hypersurfaces, say $H_1$ and $H_2$ (at least one of them bad).
\end{itemize}

The stratification at $V^{(2)}_r$ will be defined as follows: In the first case, we define a unique stratum, which is an affine line containing $\x$, say $\A^1\subset \mathbb{P}^1=H_1$. In the second case, we stratify $H_1\cup H_2$ in two strata: one affine line $\A^1$ so that $\x\in\A^1\subset H_1$  and the affine line $H_2\setminus\{\x\}$.

As we are going to blow-up along the singular locus, further quadratic transformations will be defined over $V_r^{(2)}$. This will introduce new exceptional components of the form $H_{j}=\mathbb{P}^1$. We indicate in \ref{parstr} and \ref{parinf}  how to stratify the union of new components which are bad (Definition \ref{defcolhyper}). Each stratum, will be either 
\begin{itemize}
\item[$*$] an affine line $\mathbb{A}^1$, or
\item[$*$] a point.
\end{itemize}

A zero-dimensional stratum, a point, will be called an \emph{infinitesimal stratum}. These zero-dimen\-sion\-al strata will always arise as intersection of two exceptional hypersurfaces: one good and the other bad. However, the intersection point of a bad and a good hypersurface is not necessarily a zero-dimensional stratum.
These strata will appear and treated in detail in \ref{parstr} Case C).

\end{parrafo}

\begin{parrafo}{\bf Stratification and quadratic transformations.}\label{parstr}

Once a stratification is fixed, we blow-up at a singular point, and then a new stratification will be defined. In doing that, we will follow a rule: \emph{the strict transform of an already defined stratum is a stratum}. We therefore need to establish a criterion 
to stratify points along the exceptional hypersurface, every time we blow up a point.

\vspace{0.15cm}

\noindent {\bf Case A)} In this case we assume that a quadratic transformation is defined at a point $\x$, which is bad, and is not a zero-dimensional stratum. This situation can occur within the following sub-cases:

\begin{center}
\setlength{\unitlength}{0.9mm}
\begin{picture}(140,47)
\allinethickness{1.4pt}
\put(20,5){\textcolor{darkpurple}{\line(0,1){40}}}
\put(0,25){\textcolor{darkpurple}{\line(1,0){40}}}
{\color{darkpurple}
\put(21.5,25){\makebox(0,0){$($}}
\put(21.5,25){\makebox(0,0){$($}}
\put(21.5,25){\makebox(0,0){$($}}

\put(18.6,25){\makebox(0,0){$)$}}
\put(18.6,25){\makebox(0,0){$)$}}
\put(18.6,25){\makebox(0,0){$)$}}

\put(20,25){\circle{1}}
\put(20,25){\circle{1.5}}
\put(20,25){\circle{0.4}}

}
\put(20,2){\makebox(0,0){\texttt{{\small Picture A1}}}}
\put(24,40){\makebox(0,0){$H_1$}}
\put(37.5,21.5){\makebox(0,0){$H_2$}}
\put(22.3,28.3){\makebox(0,0){$\x$}}

\put(70,5){\textcolor{darkpurple}{\line(0,1){40}}}
{\color{darkpurple}
\put(70,25){\circle{1}}
\put(70,25){\circle{1.5}}
\put(70,25){\circle{0.4}}
}
\put(70,2){\makebox(0,0){\texttt{{\small Picture A2}}}}
\put(74,40){\makebox(0,0){$H_1$}}
\put(74,25){\makebox(0,0){$\x$}}

\put(100,25){\textcolor{green}{\line(1,0){40}}}
\put(120,5){\textcolor{darkpurple}{\line(0,1){40}}}
{\color{darkpurple}
\put(120,25){\circle{1}}
\put(120,25){\circle{1.5}}
\put(120,25){\circle{0.4}}
{\color{green}
\put(121.6,25){\makebox(0,0){$($}}
\put(121.6,25){\makebox(0,0){$($}}
\put(121.6,25){\makebox(0,0){$($}}

\put(118.4,25){\makebox(0,0){$)$}}
\put(118.4,25){\makebox(0,0){$)$}}
\put(118.4,25){\makebox(0,0){$)$}}
}
}

\put(120,2){\makebox(0,0){\texttt{{\small Picture A3}}}}
\put(123,40){\makebox(0,0){$H_1$}}
\put(121.5,28.5){\makebox(0,0){$\x$}}
\put(136,21.5){\makebox(0,0){$H_2$}}

\end{picture}
\end{center}


Here, the stratification of the quadratic transformation at $\x$ will be defined in a way that will not give rise to a zero-dimensional stratum.

Along this section we agree that every time we blow-up a point, the horizontal line will denote the new exceptional component. 

Let $H_1$ denote again the strict transform of $H_1$. The new stratification along points of the exceptional component is defined as follows:

\setlength{\unitlength}{0.9mm}
\begin{center}
\begin{picture}(50,47)
\allinethickness{1.4pt}
\put(0,25){\textcolor{darkpurple}{\line(1,0){40}}}
\put(20,5){\textcolor{darkpurple}{\line(0,1){40}}}
{\color{darkpurple}
\put(21,25){\makebox(0,0){$($}}
\put(21,25){\makebox(0,0){$($}}
\put(21,25){\makebox(0,0){$($}}

\put(19.1,25){\makebox(0,0){$)$}}
\put(19.1,25){\makebox(0,0){$)$}}
\put(19.1,25){\makebox(0,0){$)$}}
}
\put(20,2){\makebox(0,0){\texttt{{\small New A-stratification after the quadratic transformation}}}}
\put(24,40){\makebox(0,0){$H_1$}}
\put(38,22){\makebox(0,0){$\mathbb{A}^1=\mathbb{P}^1\setminus\{q\}$}}
\put(22,29){\makebox(0,0){$q$}}

\end{picture}
\end{center}
\vspace{-0.2cm}

In this case, a unique $1$-dimensional stratum  $\mathbb{A}^1$ is introduced after the quadratic transformation. This new affine stratum is defined as $\mathbb{A}^1= \mathbb{P}^1\setminus\{q\}$, where here $q=H_1\cap \mathbb{P}^1$. The stratum $H_2\setminus \x$, in the pictures A1) (or in  A3)), defines a new stratum after the quadratic transformation simply by taking its strict transform. Recall that we require that the strict transform of an already defined stratum to be a stratum.

\vspace{0.15cm}

\noindent{\bf Case B)} Here, we study the case of a quadratic transformation at a good point $\x$ which is only in one bad hypersurface $H_1$, and is not a zero-dimensional stratum. This can occur in the  following sub-cases:

\setlength{\unitlength}{0.9mm}
\begin{center}
\begin{picture}(70,47)
\allinethickness{1.4pt}
\put(0,5){\textcolor{darkpurple}{\line(0,1){40}}}
{\color{green}
\put(0,25){\circle{1}}
\put(0,25){\circle{1.5}}
\put(0,25){\circle{0.4}}
}
\put(0,2){\makebox(0,0){\texttt{{\small Picture B1}}}}
\put(4,40){\makebox(0,0){$H_1$}}
\put(4,25){\makebox(0,0){$\x$}}

\put(50,5){\textcolor{darkpurple}{\line(0,1){40}}}
\put(30,25){\textcolor{green}{\line(1,0){40}}}
{\color{green}
\put(50,25){\circle{1}}
\put(50,25){\circle{1.5}}
\put(50,25){\circle{0.4}}

\put(51.5,25){\makebox(0,0){$($}}
\put(51.5,25){\makebox(0,0){$($}}
\put(51.5,25){\makebox(0,0){$($}}

\put(48.5,25){\makebox(0,0){$)$}}
\put(48.5,25){\makebox(0,0){$)$}}
\put(48.5,25){\makebox(0,0){$)$}}
}

\put(50,2){\makebox(0,0){\texttt{{\small Picture B2}}}}
\put(53,40){\makebox(0,0){$H_1$}}
\put(52.5,28.3){\makebox(0,0){$\x$}}

\end{picture}
\end{center}

\vspace{-0.2cm}

After the quadratic transformation at $\x$, the new exceptional hypersurface is good. We define the new stratification by taking the strict transform of the previous stratum defined over $H_1$: 

\setlength{\unitlength}{1mm}
\begin{center}
\begin{picture}(50,47)
\allinethickness{1.4pt}
\put(0,25){\textcolor{green}{\line(1,0){40}}}
\put(20,5){\textcolor{darkpurple}{\line(0,1){40}}}
{\color{green}
\put(21,25){\makebox(0,0){$($}}
\put(21,25){\makebox(0,0){$($}}
\put(21,25){\makebox(0,0){$($}}

\put(19.1,25){\makebox(0,0){$)$}}
\put(19.1,25){\makebox(0,0){$)$}}
\put(19.1,25){\makebox(0,0){$)$}}
}

\put(22.25,28.5){\makebox(0,0){$q$}}

\put(20,2){\makebox(0,0){\texttt{{\small New B-stratification after the quadratic transformation}}}}
\put(24,40){\makebox(0,0){$H_1$}}

\end{picture}
\end{center}

\vspace{0.15cm}

\noindent{\bf Case C)} In this case, the quadratic transformation is defined at a good point $\x$ which is an intersection of two bad lines $H_1$ and $H_2$. The point $\x$ is in the stratum defined by $H_1$, at least locally at this point. In particular, the point is not a zero-dimensional stratum.

\begin{center}
\setlength{\unitlength}{.9mm}
\begin{picture}(40,47)
\allinethickness{1.4pt}
\put(20,5){\textcolor{darkpurple}{\line(0,1){40}}}
\put(0,25){\textcolor{darkpurple}{\line(1,0){40}}}

{\color{green}
\put(20,25){\circle{1}}
\put(20,25){\circle{1.5}}
\put(20,25){\circle{0.4}}

{\color{darkpurple}
\put(21.5,25){\makebox(0,0){$($}}
\put(21.5,25){\makebox(0,0){$($}}
\put(21.5,25){\makebox(0,0){$($}}
\put(18.5,25){\makebox(0,0){$)$}}
\put(18.5,25){\makebox(0,0){$)$}}
\put(18.5,25){\makebox(0,0){$)$}}
}
}
\put(20,2){\makebox(0,0){\texttt{{\small Picture C}}}}
\put(24,40){\makebox(0,0){$H_1$}}
\put(34.5,21.5){\makebox(0,0){$H_2$}}
\put(21.5,28.5){\makebox(0,0){$\x$}}

\end{picture}
\end{center}

The new (horizontal) exceptional line is good and the stratification is defined by:

\begin{itemize}
\item The strict transform of the previous strata.
\item A new zero-dimensional stratum defined as the intersection of $H_2$ and the new exceptional good line. 
\end{itemize}
Note that in this case, depicted below, we introduce a zero-dimensional stratum $\mathcal{Q}$:

\begin{center}
\setlength{\unitlength}{.9mm}
\begin{picture}(80,45)
\allinethickness{1.4pt}
\put(0,25){\textcolor{green}{\line(1,0){80}}}
\put(20,5){\textcolor{darkpurple}{\line(0,1){40}}}

{\color{green}
\put(21,25){\makebox(0,0){$($}}
\put(21,25){\makebox(0,0){$($}}
\put(21,25){\makebox(0,0){$($}}

\put(19,25){\makebox(0,0){$)$}}
\put(19,25){\makebox(0,0){$)$}}
\put(19,25){\makebox(0,0){$)$}}
}
\put(24,40){\makebox(0,0){$H_1$}}

\put(60,5){\textcolor{darkpurple}{\line(0,1){40}}}
{\color{lightblue}
\put(60,25){\circle{3}}
\put(60,25){\circle{2.4}}
\put(60,25){\circle{1.9}}
\put(60,25){\circle{1}}
\put(60,25){\circle{1.5}}
\put(60,25){\circle{0.4}}
{\color{darkpurple}
\put(60,27){\makebox(0,0){$\smile$}}
\put(60,27){\makebox(0,0){$\smile$}}
\put(60,27){\makebox(0,0){$\smile$}}

\put(60,22.4){\makebox(0,0){$\frown$}}
\put(60,22.4){\makebox(0,0){$\frown$}}
\put(60,22.4){\makebox(0,0){$\frown$}}
}
}
\put(62,40){\makebox(0,0){$H_2$}}
\put(62.5,27){\makebox(0,0){$\mathcal{Q}$}}
\put(40,2){\makebox(0,0){\texttt{{\small New C-stratification after the quadratic transformation}}}}
\end{picture}
\end{center}

\end{parrafo}

\vspace{-0.2cm}

\begin{parrafo}\label{parinf}{\bf Quadratic transformations at an infinitesimal statum.}


\noindent{\bf Case D1)} Here we study a quadratic transformations at a bad point $\Q$, which is in addition an zero-dimensional stratum:
\vspace{-.15cm}
\begin{center}
\setlength{\unitlength}{.9mm}
\begin{picture}(40,47)
\allinethickness{1.4pt}
\put(0,25){\textcolor{green}{\line(1,0){40}}}
\put(20,5){\textcolor{darkpurple}{\line(0,1){40}}}
{\color{darkpurple}
\put(20,25){\circle{3.2}}
\put(20,25){\circle{2.5}}
\put(20,25){\circle{1.9}}
\put(20,25){\circle{1}}
\put(20,25){\circle{1.5}}
\put(20,25){\circle{0.4}}
\put(20,27){\makebox(0,0){$\smile$}}
\put(20,27){\makebox(0,0){$\smile$}}
\put(20,27){\makebox(0,0){$\smile$}}

\put(20,22.4){\makebox(0,0){$\frown$}}
\put(20,22.4){\makebox(0,0){$\frown$}}
\put(20,22.4){\makebox(0,0){$\frown$}}
}
\put(24,40){\makebox(0,0){$H_1$}}
\put(23.2,27){\makebox(0,0){$\mathcal{Q}$}}
\put(36,21.5){\makebox(0,0){$H_2$}}
\put(20,2){\makebox(0,0){\texttt{{\small Picture D1}}}}
\end{picture}
\end{center}

The new stratification after the quadratic transformation is defined by
\begin{enumerate}
\item[a)] The strict transforms of the previous strata.

\item[b)] The stratification of the new exceptional line $\mathbb{P}^1$ as the union of the  infinitesimal stratum $\Q'=\mathbb{P}^1\cap H_2$ and the $1$-dimensional stratum $\mathbb{A}^1=\mathbb{P}^1\setminus\{\Q'\}$.
\end{enumerate}

This new stratification is depicted as follows:


\begin{center}
\setlength{\unitlength}{.9mm}
\begin{picture}(80,45)
\allinethickness{1.4pt}
\put(20,5){\textcolor{green}{\line(0,1){40}}}

\put(60,5){\textcolor{darkpurple}{\line(0,1){40}}}
\put(0,25){\textcolor{darkpurple}{\line(1,0){80}}}
{\color{lightblue}
\put(20,25){\circle{3.2}}
\put(20,25){\circle{2.5}}
\put(20,25){\circle{1.9}}
\put(20,25){\circle{1}}
\put(20,25){\circle{1.5}}
\put(20,25){\circle{0.4}}
{\color{darkpurple}
\put(22.4,25){\makebox(0,0){$($}}
\put(22.4,25){\makebox(0,0){$($}}
\put(22.4,25){\makebox(0,0){$($}}

\put(17.7,25){\makebox(0,0){$)$}}
\put(17.7,25){\makebox(0,0){$)$}}
\put(17.7,25){\makebox(0,0){$)$}}

\put(60,25.8){\makebox(0,0){$\smile$}}
\put(60,25.8){\makebox(0,0){$\smile$}}
\put(60,25.8){\makebox(0,0){$\smile$}}

\put(60,23.7){\makebox(0,0){$\frown$}}
\put(60,23.7){\makebox(0,0){$\frown$}}
\put(60,23.7){\makebox(0,0){$\frown$}}
}
}
\put(62,40){\makebox(0,0){$H_1$}}
\put(22,40){\makebox(0,0){$H_2$}}
\put(21.6,29.2){\makebox(0,0){$\mathcal{Q}'$}}
\put(39,21){\makebox(0,0){$\mathbb{A}^1=\mathbb{P}^1\setminus\{\mathcal{Q}'\}$}}
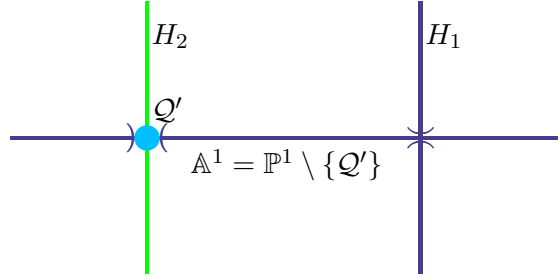
\put(40,2){\makebox(0,0){\texttt{{\small New D1-stratification after the quadratic transformation}}}}

\end{picture}
\end{center}


\noindent{\bf Case D2)} This is the case in which the center of the quadratic transformation is a good point $\Q$, which is also an infinitesimal stratum. Namely,

\begin{center}
\setlength{\unitlength}{.9mm}
\begin{picture}(40,47)
\allinethickness{1.4pt}
\put(0,25){\textcolor{green}{\line(1,0){40}}}
\put(20,5){\textcolor{darkpurple}{\line(0,1){40}}}
{\color{green}
\put(20,25){\circle{3.2}}
\put(20,25){\circle{2.5}}
\put(20,25){\circle{1.9}}
\put(20,25){\circle{1}}
\put(20,25){\circle{1.5}}
\put(20,25){\circle{0.4}}
{\color{darkpurple}
\put(20,27.1){\makebox(0,0){$\smile$}}
\put(20,27.1){\makebox(0,0){$\smile$}}
\put(20,27.1){\makebox(0,0){$\smile$}}

\put(20,22.3){\makebox(0,0){$\frown$}}
\put(20,22.3){\makebox(0,0){$\frown$}}
\put(20,22.3){\makebox(0,0){$\frown$}}
}
}
\put(24,40){\makebox(0,0){$H_1$}}
\put(36,21.5){\makebox(0,0){$H_2$}}
\put(23.2,27){\makebox(0,0){$\mathcal{Q}$}}
\put(20,2){\makebox(0,0){\texttt{{\small Picture D2}}}}
\end{picture}
\end{center}

\vspace{0.2cm}

The new stratification is represented by

\begin{center}
\setlength{\unitlength}{.9mm}
\begin{picture}(80,47)
\allinethickness{1.4pt}
\put(20,5){\textcolor{green}{\line(0,1){40}}}

\put(60,5){\textcolor{darkpurple}{\line(0,1){40}}}
\put(0,25){\textcolor{green}{\line(1,0){80}}}

{\color{lightblue}
\put(60,25){\circle{3.2}}
\put(60,25){\circle{2.5}}
\put(60,25){\circle{1.9}}
\put(60,25){\circle{1}}
\put(60,25){\circle{1.5}}
\put(60,25){\circle{0.4}}
{\color{darkpurple}
\put(60,27.1){\makebox(0,0){$\smile$}}
\put(60,27.1){\makebox(0,0){$\smile$}}
\put(60,27.1){\makebox(0,0){$\smile$}}

\put(60,22.3){\makebox(0,0){$\frown$}}
\put(60,22.3){\makebox(0,0){$\frown$}}
\put(60,22.3){\makebox(0,0){$\frown$}}
}
}
\put(63,40){\makebox(0,0){$H_1$}}
\put(22,40){\makebox(0,0){$H_2$}}
\put(63,27.4){\makebox(0,0){$\mathcal{Q}'$}}
{\color{green}
\put(20,25){\makebox(0,0){$($}}
\put(20,25){\makebox(0,0){$($}}
\put(20,25){\makebox(0,0){$($}}

\put(18,25){\makebox(0,0){$)$}}
\put(18,25){\makebox(0,0){$)$}}
\put(18,25){\makebox(0,0){$)$}}
}
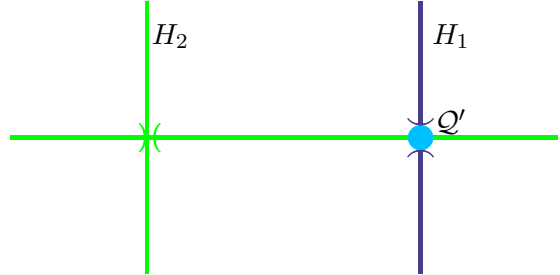
\put(40,2){\makebox(0,0){\texttt{{\small New D2-stratification after the quadratic transformation}}}}
\end{picture}
\end{center}
So only a new zero-dimensional stratum $\mathcal{Q}'$ is introduced at the bad locus.

\end{parrafo}

\begin{parrafo}{\bf Definition of the local data and local invariants}.\label{par:definv}

\vspace{0.1cm}

Once the stratification has been fixed, notions of {\em local data} and {\em local invariants} will be introduced at every isolated point $x\in\Sing(\G_{r'})$ (or equivalently, to $\x\in\beta_{r'}(\Sing(\G_{r'}))$. Here, the local data will assign to each $1$-dimensional stratum $\A^1=\Spec(k[y])$ a polynomial $g(y)$ in $k[y]$,  and to each zero-dimensional stratum $\mathcal{Q}$ an element of $\calo_{V^{(2)}_{r'},\mathcal{Q}}$. Finally, local invariants will be defined in terms of these local data.

\vspace{0.3cm}

\noindent $\bullet$ {\bf Local data in case $\x \in {\mathbb A}^1$:} Here we take an isolated point $\x\in\beta_{r'}(\Sing(\G_{r'}))$ in a $1$-dimensional stratum  ${\mathbb A}^1$ (included in the union of the bad hypersurfaces, and in particular located in a bad hypersurface, say $H_1$). 
Let $p\P=p\P(\beta_{r'},z,f_{p^e})$ be a well-adapted $p$-presentation at $\x$, where 
\begin{equation}\label{eqinvafpe}
f_{p^e}(z)=z^{p^e}+a_{1}z^{p^e-1}+\dots+a_{p^e}\in\calo_{V^{(2)}_{r'},\x}[z].
\end{equation}
 Fix a regular system of parameters $\{x,y\}$ at $\calo_{V_{r'}^{(2)},\x}$, so that $\{x=0\}$ defines $\A^1$ locally at $\x$.

As $H_1$ is a bad hypersurface, $\Hord^{(2)}(\G_{r'})(\xi_{H_1})=\frac{\nu_{\xi_{H_1}}(a_{p^e})}{p^e}<\ord((\R_{\G,\beta})_{r'})(\xi_{H_1})$, where $\xi_{H_1}$ is the generic point of $H_1$. Note that, in particular, $p^e\cdot\Hord^{(2)}(\G_{r'})(\xi_{H_1})\in\mathbb{Z}_{\geq0}$. Let us denote this integer by $\ell=p^e\cdot \Hord^{(2)}(\G_{r'})(\xi_{H_1})$. There is a factorization of the form
\begin{equation}\label{lec}
a_{p^e}=x^\ell(g(y)+x\Omega(x,y)),
\end{equation}
where $\Omega(x,y)\in\calo_{V_{r'}^{(2)},\x}$, and the exponent $\ell<p^e$ by (\ref{arbol}) . Note here that $H_1$ is an exceptional hypersurface introduced by previous quadratic transformation, it can be checked that local coordinates $x$ and $y$ can be chosen so that $g(y)$ is indeed a polynomial in $k[y]$.

In this case, the \emph{local data} at $\x$ will be defined as the pair $\big(\frac{\ell}{p^e},g(y)W^{p^e}\big).$ 

\begin{definition}\label{papel} Fix, with the setting as above, a point $\x\in\A^1$, and the local data $\big(\frac{\ell}{p^e},g(y)W^{p^e}\big).$ The \emph{order of $g$ at $\x$}, say $\ord_\x(g)$, is defined as follows: 
\begin{itemize}
\item If $\ell\not=0$, $\ord_{\x}(g)$ is the usual order of $g$ at $\calo_{\A^1,{\x}}$.
\item If $\ell=0$, $\ord_{\x}(g)$ is the smallest power of $y$, which appears in the Taylor expansion of $g(y)$ at the point, that is not a $p^e$-th power.
\end{itemize}
\end{definition}

\vspace{0.2cm}

\noindent$\bullet$ {\bf Local data in case $x =\mathcal{Q}$ is an infinitesimal statum:} Set local coordinates so that  $(\R_{\G,\beta})_{r'}=x^ay^bW^s$, assuming that $H_1=\{x=0\}$ denotes the good hypersurface through the point. Now, $\Hord^{(2)}(\G_{r'})(\xi_{H_1})=\ord((\R_{\G,\beta})_r)(\xi_{H_1})=\frac{a}{s}$. Define the \emph{local data} as the pair 
$\big(\frac{a}{s},g(y)W^s\big),$ where $g(y)=y^b$.

\begin{definition}({\em Local invariants}).
Fix a point $\x\in\beta_{r'}(\Sing(\G_{r'}))$ with local data, say $(\frac{m}{n},g(y)W^t)$, for a suitable integer $t$. The \emph{local invariant} assign to $\x$ is defined as:
\begin{itemize}
\item If $\x\in\A^1$, then $\inv(\x)=\frac{\ord_\x(g)}{p^e}$ (in this case $t=p^e$).
\item If $\x=\mathcal{Q}$, then $\inv(\x)=\frac{\nu_\x(g)}{s}$ (here $t=s$, and $g(y)W^t=y^bW^s$).
\end{itemize}
\end{definition}
\end{parrafo}

\vspace{0.15cm}

\begin{parrafo}{\bf Invariants and transformations}.\label{par:invtrans}

We now study the behavior of the previous invariants under quadratic transformations, taking into account the distinction in the cases presented in  \ref{parstr}.

\vspace{0.3cm}

\noindent{\bf Case A)} 
Local coordinates $\{x,y\}$ are chosen locally at $\calo_{V_{r'}^{(2)},\x}$ so that locally $\A^1=\{x=0\}$. 
In this case, as $\x$ is assumed to be bad, and a well-adapted $p$-presentation can be chosen so as
$$\Hord^{(2)}(\G_{r'})(x)=\frac{\nu_\x(a_{p^e})}{p^e}<\ord((\R_{\G,\beta})_{r'})(\x).$$
In addition, the initial form of $a_{p^e}$ in $\Gr_\x(\calo_{V_{r'}^{(2)}})$ is not a $p^e$-th power.

The objective is to define local invariants after the quadratic transformation at $\x$. This leads to:
\begin{enumerate} 
\item The definition of local data and invariants at the strict transform of $H_1$.
\item The definition of local data and invariants at points in the new stratum $\A^1$ (included in the exceptional component).
\end{enumerate}

\vspace{0.3cm}

\noindent{\bf (1)} {\it Local invariants at the strict transform of $\mathbb{A}^1\subset H_1$}.
\vspace{0.1cm}

Local factorization at $\x$  is $a_{p^e}=x^\ell(g(y)+x\Omega(x,y))$. At $\x_1$, the origin of the $U_y$-chart, (with local coordinates $x_1=\frac{x}{y}$, $y_1=y$), the factorization is given by
$$a_{p^e}^{(1)}=x_1^\ell y_1^{\ell-p^e}(g(y_1)+x_1y_1\Omega').$$
In this case, the new local data, say $\big(\frac{\ell_1}{p^e},g_1(y_1)W^{p^e}\big)$, is defined by setting $g_1(y_1)=y_1^{\ell-p^e}\cdot g(y_1)$ and $\ell_1=\ell$.

As we assume that the tight monomial algebra $\m_{r'} W^s$ is reduced (see \ref{par:mred}), then $\ell<p^e$, and hence  
\begin{equation}\label{eqinvH1drop}\ord_{\x_1}(g_1)<\ord_\x(g).
\end{equation}

\vspace{0.3cm}

\noindent{\bf (2)} {\it Invariants along the new stratum $\A^1$}.
\vspace{0.1cm}

Set $\mathbb{A}^1=\mathbb{P}^1\setminus\{\x_1\}=\Spec(k[Y])$ to be the 1-dimensional stratum defined as in Case A) in \ref{parstr}. We first assign a polynomial $g_1(Y)$ to $\A^1$.

Locally at $\x$, the $1$-dimensional stratum  containing $\x$ was defined by $\{x=0\}$. The new stratum along the exceptional hypersurface, the affine line $\A^1$, is the intersection of the exceptional hypersurface with the open chart $U_{x_1}$ (with coordinates  $x_1=x$ and $y_1=\frac{y}{x}$).

Set $\Gr_\x(\calo_{V^{(2)}_{r'}})=k'[X,Y]$ with $X=\In_\x(x)$ and $Y=\In_\x(y)$, and set $\In_\x(a_{p^e})=\sum_{i+j=d}\lambda_{i,j}X^iY^j,$
the initial form of $a_{p^e}$ at $\mathbf{x}$, where $d=\nu_\x(a_{p^e})$.
Finally, define $\widetilde{g}(Y)=\In_\x(a_{p^e})|_{X=1}$.

Fix a point in $\A^1$, and after suitable change of coordinates, say $y_1=Y+\alpha$, assume that $y_1$ vanishes at such point. Now set $g_1(y_1)=\widetilde{g}(Y+\alpha)$. Then, the local invariant is 
$(\frac{\ell_1}{p^e},g_1(y_1))$,
with $\frac{\ell_1}{p^e}=\Hord^{(2)}(\G_{r'+1})(\xi_{H'})$, where $\xi_{H'}$ is the generic point of the new exceptional hypersurface, say $H'$.
The previous change of variables does not affect the degree of the polynomial $g_1$. Namely $\deg(\widetilde{g}(Y))=\deg(g_1(y_1))$.

\begin{lemma}\label{leminvAdes} {\rm(}{\bf Abhyankar's trick}{\rm )}. Fix $x\in\Sing(\G_{r'})$ (or equivalently $\x\in\beta_{r'}(\Sing(\G_{r'}))$) and assume that the setting is as above. For any point $\mathbf{x}'\in\mathbb{A}^1\subset H'$,
\begin{equation}\label{eqinvAdes}
\ord_{\x'}(g_1)\leq \ord_\x(g).
\end{equation}
\end{lemma}

\begin{proof} 
Consider $a_{p^e}$ as a formal power serie in the variables $x$ and $y$,  (\ref{lec}) indicates that it is expressed as a sum of monomials of the form $x^t y^r$, with $t\geq \ell$. Note, in addition, that a monomial of the form $x^\ell y^M$, where $M=\ord_\x(g(y))$, appears in such formal expression.

This leads to the following conclusions:
\begin{enumerate}
\item $\In_\x(a_{p^e})=\sum_{i\geq \ell}\lambda_{i,d-i} X^iY^{d-i}.$
\item $d \leq \ell+M= \deg (x^\ell y^M)$.
\end{enumerate}
From where it is inferred that if $\lambda_{i,d-i} \neq 0$, then 
$ d-i\leq M$. In particular,  $\sum_{i\geq \ell}\lambda_{i,d-i} Y^{d-i}$ is a polynomial of degree $M_1$ with  $ M_1 \leq \ord_\x(g)=M$. We {\bf claim} that $\ord_{\x'}(g_1)\leq M_1$, and this would ensure that
 the inequality (\ref{eqinvAdes}). The proof of this claim will be addressed in \ref{tulsa}, it will make use of the following Lemma.
 \end{proof}
 
 \begin{lemma}\label{estrella}
Fix a polynomial $g(y)\in k[y]$. Then at any closed point $\x\in\mathbb{A}^1_k=\Spec(k[y])$:
\begin{enumerate}
\item $\ord_\x(g)\leq \deg(g)$, where $\ord_\x(g)$ is the usual order of $g$ at $\calo_{\mathbb{A}^1_k,\x}$.
\item If $g(y)\not\in k[y^{p^e}]$, then $p^e$-$\ord_\x(g)\leq \deg(g)$, where $p^e$-$\ord_\x(g)$ denotes the smallest power of $y$, that is not a $p^e$-th power, which appears in the Taylor expansion of $g(y)$ at the point.
\end{enumerate}
 \end{lemma}
 
\begin{proof} Let $M$ denote the degree of $g(y)$.

(1) \ Fix a change of variables $y_1=y-\alpha$ so that $y_1$ vanishes at $\x$, then $g(y)=g_1(y_1)$ is also a polynomial of degree $M$ in $y_1$. Hence, $\ord_\x(g)\leq \deg(g)$.

(2) Let $M'\leq M$ be the biggest integer so that $M'\not\equiv 0\mod p^e$ and that $y^{M'}$ appears in $g(y)$. 
Consider the Taylor expansion of $g(y)$ as in (\ref{lapiz}) applied here for $S[Z]=k[y]$. As $\Delta^{M'}(g)\in k\setminus\{0\}$, then the term $y_1^{M'}$ appears at the Taylor development at a fixed point $\x$. Here $y_1=y-\alpha$ is defined so as to vanish at $\x$.
\end{proof}

\begin{parrafo}\label{tulsa} We address now the proof of the claim stated in Lemma \ref{leminvAdes}. Fix notation as in this lemma, where $\x\in\Sing(\G_{r'})$ and $H'$ is the exceptional hypersurface introduced by blowing-up $\x$.
The stratum $\mathbb{A}^1\subset H'$ is defined by $\A^1=\Spec(k[y_1])$, where $y_1=\frac{y}{x}$. We consider $g_1(y_1)\in k[y_1]$ to be naturally identified with $\In_\x(a_{p^e})|_{x=1}=\sum_{i\geq \ell}\lambda_{i,d-i}y_1^{d-i}$, so $g_1(y_1)$ is a polynomial of degree $M_1\leq M$.
Recall now the notion of invariant attached to a singular point $\x'\in\mathbb{A}^1$ in Definition \ref{papel}.

Set as before $d=\nu_\x(a_{p^e})$, locally at $\x$.  We distinguish two cases:

(a) If $d\not\equiv 0\mod p^e$, then the local data at $\x'$ is $(\ell_1,g_1(y_1))$ with $\ell_1\not=0$, and hence $\ord_{\x'}(g_1)$  is the usual order at $\calo_{\mathbb{A}^1,\x'}$.  Inequality (\ref{eqinvAdes}) follows from Lemma \ref{estrella} (1).

(b) If $d\equiv 0 \mod p^e$, then the local data at $\x'$ is $(\ell_1,g_1(y_1))$ with $\ell_1=0$. Here $\ord_{x'}(g_1)$ is provided by $p^e$-$\ord_{\x'}(g_1)$ as in Lemma \ref{estrella} (2), and the condition $g_1(y_1)\in k[y_1]$ is guaranteed by the fact that $\In_{\x}(a_{p^e})$ is not a $p^e$-power.
\end{parrafo}

\begin{lemma}\label{lem:notseq}
If there is a sequence of quadratic transformations at points $\x_0,\x_1\dots$, where each $\x_i$ maps to $\x_{i-1}$ and such that:
\begin{enumerate}
\item Case {\rm A)} occurs at each point $\x_i$, and
\item $\ord_{\x_i}(g_{i})=\ord_{\x_{i-1}}(g_{i-1})$, 
\end{enumerate}
then the sequence must be finite.
\end{lemma}

\begin{proof}

The existence of an infinite sequence with the previous conditions would contradict the assumption that $\beta_{r'}(\Sing(\G_{r'}))$ has no $1$-codimensional component. In fact, assume that such a sequence does exist. We claim that at the completion at the closed point $\x_0$, there is a smooth curve whose successive strict transforms passes through the sequence $\x_i$. This follows from the assumption that (\ref{eqinvH1drop}) does not occur. This ensures that there is a smooth curve passing through these points, and a smooth curve with this property would be a $1$-dimensional component of the singular locus. This contradicts the hypothesis. 
\end{proof}

We assume that  $\beta_{r'}(\Sing(\G_{r'}))$ is a finite set of closed points. Every time we fix one such point, there is a unique procedure of quadratic transformations over it. In fact, after applying a quadratic transformation, we reduce to the case $\beta_{r'+1}(\Sing(\G_{r'+1}))$ is a finite set of closed points by blowing-up the new exceptional component. Lemma \ref{lem:notseq} ensures that over $\x_0\in \beta_{r'}(\Sing(\G_{r'}))$ only finitely many singular points  $\x_0,\x_1\dots$ can arise, where each $\x_i$ maps to $\x_{i-1}$, and with the condition that Case A) is preserved and equality holds in (\ref{eqinvAdes}).

\begin{lemma}\label{lemfinite}
Fix $\x_0\in \beta_{r'}(\Sing(\G_{r'}))$, there is a uniform bound for the length of all possible sequences of quadratic transformations over $\x_0$ in the setting of Lemma \ref{lem:notseq}.
\end{lemma}

\begin{proof}
Fix a point $\x\in\beta_{r'}(\Sing(\G_{r'}))$ which we may assume to be isolated and within case A). We claim that after a finite sequence of quadratic transformations over this point, the invariants drops at {\em any exceptional point} $\x_i$ mapping to $\x$.

To this end we first show that after finitely many blow ups any singular point mapping to $\x$, for which condition A) holds and equality occurs at (\ref{eqinvAdes}),  must be in case A2).

To check this claim note first that under the assumption of the equality in (\ref{eqinvAdes}), case A2) is stable. Namely, if case A2) holds at a point, and case A) holds at a point after the quadratic transformation, and equality holds at 
(\ref{eqinvAdes}), then this exceptional point must also be within case A2). On the other hand, assuming that $\x$ is in condition A1) or A3), for which $\x$ is contained in the hypersurface $H_2$, Lemma \ref{lem:notseq} ensures that after finitely many quadratic transformations, every point within case A) and for which equality holds, must be in case A2). In fact, otherwise the exceptional hypersurface $H_2$ would be a component of $\beta_{r'}(\Sing(\G_{r'}))$, in contradiction with our hypothesis. This last assertion follows using same arguments as before.

The previous finite sequence of quadratic transformations over the point $\x$, constructed so as to be in case A2), introduces finitely many new exceptional components, say $H_{n_1},\dots,H_{n_t}$. Let the elimination algebra be of the form
$(\R_{\G,\beta})_{n_t}=I(H_1)^{\alpha_1}\dots I(H_{n_t})^{\alpha_{n_t}}W^s,$
for some integers $\alpha_i\geq 0$. 
We claim now that after at most $\alpha_1+\dots+\alpha_{n_t}$ quadratic transformations, the inequality (\ref{eqinvAdes}) will be strict at any point mapping to $\x$ which fulfills condition A). 

To check this, note first that locally at any point within condition A2),  there is a regular system of parameters $\{x,y\}$ so that $(\R_{\G,\beta})_{n_t}=x^aW^s$ with $a=\alpha_i$ for some $i\in\{1,\dots,n_t\}$. Finally, note that a quadratic transformation at each point introduces a new hypersurface, that any exceptional singular point is in case A2), and that $(\R_{\G,\beta})_{n_t+1}=x_1^{a-s}W^s$. This can occur only finitely many times.
The claim follows now from the inclusion $\beta_j(\Sing(\G_j))\subset\Sing((\R_{\G,\beta})_j)$.
\end{proof}

\vspace{0.25cm}

\noindent{\bf Case B)}
In this case, the stratification is defined by the strict transform of the previous stratum. So attention should be drawn only at the unique point $q$ of the strict transform of $H_1$.

This parallels the situation of case A) (1) and the invariant strictly drops as in (\ref{eqinvH1drop}).

\vspace{0.25cm}

\noindent{\bf Case C)} 
If a point $\x$ is within case C), then $\x$ is an intersection of two bad exceptional hypersurfaces $H_1$ and $H_2$. Moreover, the point $\x$ is good and belongs to a $1$-dimensional stratum $\mathbb{A}^1$ included in $H_1$.

Therefore a quadratic transformation at such point $\x$, introduces a good exceptional hypersurface $\mathbb{P}^1$. Locally over $\x$ the new stratification is defined by:
\begin{itemize}
\item The strict transform of the previous strata.
\item A zero-dimensional stratum $\mathcal{Q}=\mathbb{P}^1\cap H_2$.
\end{itemize}

The invariants at the strict transform of $H_1$ are to be dealt with exactly as in \ref{par:invtrans} Case A)  (1). We therefore restrict attention to the data and invariants to be defined at $\mathcal{Q}$.

Let us fix locally at $\x$ coordinates  $x, y$ so that $H_1=\{x=0\}$ and $H_2=\{y=0\}$. Assume that a local presentation is given so that $a_{p^e}$ is as in (\ref{lec}). Therefore, the local invariant at $\x$ is $\big(\frac{\ell}{p^e},g(y)W^{p^e}\big)$. Set  $(\R_{\G,\beta})_{r'}=x^ay^bW^s$. As we assume that the point $\x$ is good, then $\frac{a+b}{s}\leq \frac{\nu_\x(a_{p^e})}{p^e}.$

The point $\mathcal{Q}$ is the origin at the $U_x$-chart (with coordinates $x_1=x$, $y_1=\frac{y}{x}$). Consider the quadratic transformation at $\x$. The new  exceptional line is good, and hence the exponents of the tight monomial algebra and the elimination algebra along this hypersurface coincide. So, after reduction, we may assume that at $\mathcal{Q}$, $(\R_{\G,\beta})_{r'+1}=x_1^{a+b-sm}y_1^{b}W^s,$
for a suitable integer $m\geq 0$ so that $h_{r'+1}=a+b-sm<s$ (here $h_{r'+1}$ is also the exponent in the tight monomial algebra of the new exceptional hypersurface).

Set $g_1(y_1)W^s=y_1^bW^s$. According to \ref{par:definv}, case $x=\mathcal{Q}$, the local data we assign to $\mathcal{Q}$ is $\Big(\frac{a+b-sm}{s}, y_1^bW^s\Big).$

\begin{lemma} Assume that the conditions in the previous setting hold. Then $\frac{b}{s}=\inv(\mathcal{Q})<\inv(\x).$
\end{lemma}

\begin{proof}
Set $M=\ord_\x(g)$ and recall that $\inv(\x)=\frac{M}{p^e}$. Then
$$\frac{\ell+M}{p^e}\geq \frac{\nu_\x(a_{p^e})}{p^e}\geq \frac{a+b}{s}>\frac{\ell}{p^e}+\frac{b}{s}.$$
The first inequality follows from the fact that $x^\ell y^M$ is a monomial that appears in the formal expansion of $a_{p^e}$. The second inequality is due to the fact that $\x$ is a good point. Finally, the last inequality is a consequence of the fact that $H_1$ is bad and hence $\frac{a}{s}>\frac{\ell}{p^e}$. These inequalities imply that $\frac{M}{p^e}>\frac{b}{s}$.
\end{proof}

\vspace{0.2cm}

\noindent{\bf Case D1)} In this case, we blow-up a bad point $\Q$, which is an infinitesimal stratum, and hence $\Q$ is the intersection of a bad hypersurface $H_1$ and a good hypersurface $H_2$. As $\Q$ is bad, the quadratic transformation at $\Q$ will introduce a exceptional line $\mathbb{P}^1$, which is bad. It will give rise to two strata: and infinitesimal stratum $\mathcal{Q}'$ and an affine line $\A^1=\mathbb{P}^1\setminus\{\mathcal{Q}'\}$.
Fix local coordinates $\{x,y\}$ at $\mathcal{Q}$ so that $\{x=0\}$ defines the bad line $H_1$ and $\{y=0\}$ defines the good line $H_2$. Set $(\R_{\G,\beta})_{r'}=x^ay^bW^s$.

\begin{lemma}
Fix a point $\x'\in\A^1$. The invariant strictly drops, i.e., $\inv(\x')<\inv(\mathcal{Q}).$
\end{lemma}

\begin{proof}
At $\mathcal{Q}$ the second coordinate of the local data is given by $g(x)W^s=x^aW^s$. So the local invariant at $\mathcal{Q}$ is $\inv(\Q)=\frac{\nu_\Q(g)}{s}=\frac{a}{s}.$

As we assume that $\Q$ is a bad point, $\frac{d}{p^e}:=\frac{\nu_\mathcal{Q}(a_{p^e})}{p^e}<\ord((\R_{\G,\beta})_{r'})(\Q)=\frac{a+b}{s},$ where now
$a_{p^e}$ is as in (\ref{eqinvafpe}). On the other hand, $y=0$ defines the good line, so we claim that $\frac{b}{s}<1$. In fact, good hypersurfaces are, by definition, those for which the corresponding exponents at the elimination algebra, and at the tight monomial algebra, coincide. Since we assume that the tight monomial algebra is reduced (see \ref{par:mred}), $\frac{b}{s}<1$.

Let $\In_{\mathcal{Q}}(a_{p^e})=\sum_{i+j=d}\lambda_{i,j} X^iY^j$ denote the initial form of $a_{p^e}$ at $\Q$. Note that $\frac{j}{p^e}\geq\frac{b}{s}$.
From the previous inequalities, we obtain 
$$\frac{a+b}{s}>\frac{d}{p^e}=\frac{i+j}{p^e}\geq\frac{i}{p^e}+\frac{b}{s},$$
and hence, $\frac{a}{s}>\frac{i}{p^e}$.

At the $U_y$-chart (with coordinates $y_1=y$, $x_1=\frac{x}{y}$),  
$a_{p^e}^{(1)}=y_1^{d-p^e}(g_1(x_1)+y_1\Omega'),$
where $g_1(x_1)$ is obtained by the global polynomial $\In_\Q(a_{p^e})|_{Y=1}$. So the previous discussion shows that if $\lambda_{i,j}\not=0$, then  $\frac{a}{s}>\frac{i}{p^e}$. This, in turn, suffices to check that $\inv(\x')=\frac{\ord_{\x'}(g_1)}{p^e}<\frac{a}{s}=\inv(\Q).$
\end{proof}

Now, we study the invariant at $\Q'$. Note that local coordinates at  $\mathcal{Q}'$ are given by $x_1=x$, $y_1=\frac{y}{x}$. So the elimination algebra is $(\R_{\G,\beta})_{r'+1}=x_1^{a+b-s}y^b_1W^s$ and the strict transform of the good hypersurface $H_2$ is given by $\{ y_1=0\}$. Therefore, the second coordinate of the local data at $\Q'$ is  $g_{1}(x)W^s=x_1^{a+b-s}W^s$.
\begin{lemma} With the previous setting, $\inv(\Q')<\inv(\Q).$
\end{lemma}

\begin{proof} Recall that the second coordinate of the local data at $\Q$ is $g(x)W^s=x^aW^s$. 
By definition $\displaystyle\inv(\Q)=\frac{\nu_\Q(g)}{s}=\frac{a}{s}$ and $\displaystyle\inv(\Q')=\frac{\nu_{\Q'}(g_1)}{s}=\frac{a+b-s}{s}$. Since $H_2$ is a good hypersurface, then $\frac{b}{s}<1$, from which the strict inequality is clear.
\end{proof}

\vspace{0.3cm}

\noindent{\bf Case D2)} This is, as in D1), the case of a quadratic transformation at a point $\Q$ which is an infinitesimal stratum. It is the intersection of a bad hypersurface $H_1$ and a good hypersurface $H_2$. We assume now, in addition, that $\Q$ is a good point. Since the new exceptional hypersurface is good, there is a unique stratum $\mathcal{Q}'$ which will be infinitesimal. 

Fix local coordinates $\{x,y\}$ at $\Q$, so that  $\{x=0\}$ defines the bad line $H_1$ and $\{y=0\}$ defines the good line $H_2$. Set $(\R_{\G,\beta})_{r'}=x^ay^bW^s$. Recall that the second coordinate of the local data is  $g(x)W^s=x^aW^s$. 

Note that $x_1=\frac{x}{y},y_1=y$ are local coordinates at $\mathcal{Q}'$. The elimination algebra is $(\R_{\G,\beta})_{r'+1}=x_1^{a}y_1^{a+b-s}W^s$. The new good exceptional hypersurface is defined by $y_1=0$, so the second coordinate of the local data is $g_1(x_1)W^s=x_1^aW^s$.

\begin{remark} If case D2) holds, then
\begin{itemize}
\item $\nu_{\Q'}(g_1)=\nu_{\Q}(g)$ and hence $\inv(\Q')=\inv(\Q)$. 
\item On the other hand, 
$\frac{a+b-s}{s}=\frac{a}{s}+\Big(\frac{b}{s}-1\Big)<\frac{a}{s}.$
\end{itemize}
We conclude that case D2) cannot occur in a successive manner more than finitely many times. So, in particular, after finitely many quadratic transformations at infinitesimal strata, case D1) holds.
\end{remark}
\end{parrafo}

\newpage

In the following Table we indicate, in a synthetic manner, why resolution is achieved.

\centering
\begin{tabular}[t]{| p{5cm}| p{5cm} | p{6cm} |} 

\hline

\centering{\small{\bf
Initial stratification}} & 
\qquad\quad{\small{\bf After blow-up}} & 
\qquad\qquad\quad{\small{\bf Invariants}}
\\ 
\hline

\vspace{0.05cm}
{\bf CASE A)}
\vspace{-0.1cm}
\setlength{\unitlength}{0.9mm}
\begin{center}
\begin{picture}(10,40)
\allinethickness{1.4pt}

\put(5,2){\textcolor{darkpurple}{\line(0,1){40}}}
{\color{darkpurple}
\put(5,23){\circle{1}}
\put(5,23){\circle{1.5}}
\put(5,23){\circle{0.4}}
}
\put(9,39){\makebox(0,0){$H_1$}}
\put(8.6,23){\makebox(0,0){$\x$}}
\end{picture}
\end{center}
\vspace{-0.9cm}

&   

\vspace{0.45cm}

\

\vspace{-0.15cm}
\begin{center}
\setlength{\unitlength}{.9mm}
\begin{picture}(45,40)
\allinethickness{1.4pt}
\put(0,22){\textcolor{darkpurple}{\line(1,0){40}}}
\put(20,2){\textcolor{darkpurple}{\line(0,1){40}}}
{\color{darkpurple}
\put(21,22){\makebox(0,0){$($}}
\put(21,22){\makebox(0,0){$($}}
\put(21,22){\makebox(0,0){$($}}

\put(19.1,22){\makebox(0,0){$)$}}
\put(19.1,22){\makebox(0,0){$)$}}
\put(19.1,22){\makebox(0,0){$)$}}
}
\put(24,38){\makebox(0,0){$H_1$}}
\put(35,18){\makebox(0,0){$\mathbb{A}^1=\mathbb{P}^1\setminus\{\x_1\}$}}
\put(23,26){\makebox(0,0){$\x_1$}}
\end{picture}
\end{center}
\vspace{-0.9cm}
&

\vspace{0.3cm}
$\bullet$ After finitely many quadratic transformations, either the invariant improves (strictly drops), or it remains equal. In this last case, the singular point must be either in case  B) or C). (Lemma \ref{leminvAdes}, Lemma \ref{lem:notseq}, and Lemma \ref{lemfinite}).
 \\ 

\hline

\vspace{0.1cm}

{\bf CASE B)}
\begin{center}
\setlength{\unitlength}{0.9mm}
\begin{picture}(10,40)
\allinethickness{1.4pt}
\put(5,5){\textcolor{darkpurple}{\line(0,1){40}}}
{\color{green}
\put(5,25){\circle{1}}
\put(5,25){\circle{1.5}}
\put(5,25){\circle{0.4}}
}

\put(9,40){\makebox(0,0){$H_1$}}
\put(8.6,25){\makebox(0,0){$\x$}}
\end{picture}
\end{center}
\vspace{-1.2cm}
&

\vspace{0.52cm}

\setlength{\unitlength}{.9mm}
\begin{center}
\begin{picture}(40,40)
\allinethickness{1.4pt}
\put(0,25){\textcolor{green}{\line(1,0){40}}}
\put(20,5){\textcolor{darkpurple}{\line(0,1){40}}}
{\color{green}
\put(21,25){\makebox(0,0){$($}}
\put(21,25){\makebox(0,0){$($}}
\put(21,25){\makebox(0,0){$($}}

\put(19.1,25){\makebox(0,0){$)$}}
\put(19.1,25){\makebox(0,0){$)$}}
\put(19.1,25){\makebox(0,0){$)$}}
}
\put(24,40){\makebox(0,0){$H_1$}}
\put(23.3,27.8){\makebox(0,0){$\x_1$}}
\end{picture}
\end{center}
\vspace{-1.2cm}

&

\vspace{1.1cm}

$\bullet$ This case adds no new stratum.

$\bullet$ The invariant strictly drops, i.e., 
\quad $\inv(\x)<\inv(\x_1).$
\\

\hline

\vspace{0.1cm}

{\bf CASE C)}
\vspace{-0.3cm}
\begin{center}
\setlength{\unitlength}{.9mm}
\begin{picture}(40,40)
\allinethickness{1.4pt}
\put(20,0){\textcolor{darkpurple}{\line(0,1){40}}}
\put(0,20){\textcolor{darkpurple}{\line(1,0){40}}}

{\color{green}
\put(20,20){\circle{1}}
\put(20,20){\circle{1.5}}
\put(20,20){\circle{0.4}}

{\color{darkpurple}
\put(21.5,20){\makebox(0,0){$($}}
\put(21.5,20){\makebox(0,0){$($}}
\put(21.5,20){\makebox(0,0){$($}}
\put(18.5,20){\makebox(0,0){$)$}}
\put(18.5,20){\makebox(0,0){$)$}}
\put(18.5,20){\makebox(0,0){$)$}}
}
}
\put(23,37){\makebox(0,0){$H_1$}}
\put(35,17){\makebox(0,0){$H_2$}}
\put(21,23.2){\makebox(0,0){$\x$}}

\end{picture}
\end{center}
\vspace{-0.7cm}

 & 
 \vspace{0.2cm}

 \begin{center}
\setlength{\unitlength}{.9mm}
\begin{picture}(80,40)
\allinethickness{1.4pt}
\put(0,20){\textcolor{green}{\line(1,0){50}}}
\put(13,0){\textcolor{darkpurple}{\line(0,1){40}}}

{\color{green}
\put(14,20){\makebox(0,0){$($}}
\put(14,20){\makebox(0,0){$($}}
\put(14,20){\makebox(0,0){$($}}

\put(12,20){\makebox(0,0){$)$}}
\put(12,20){\makebox(0,0){$)$}}
\put(12,20){\makebox(0,0){$)$}}
}
\put(17,35){\makebox(0,0){$H_1$}}

\put(37,0){\textcolor{darkpurple}{\line(0,1){40}}}
{\color{lightblue}
\put(37,20){\circle{3}}
\put(37,20){\circle{2.4}}
\put(37,20){\circle{1.9}}
\put(37,20){\circle{1}}
\put(37,20){\circle{1.5}}
\put(37,20){\circle{0.4}}
{\color{darkpurple}
\put(37,22){\makebox(0,0){$\smile$}}
\put(37,22){\makebox(0,0){$\smile$}}
\put(37,22){\makebox(0,0){$\smile$}}

\put(37,17.4){\makebox(0,0){$\frown$}}
\put(37,17.4){\makebox(0,0){$\frown$}}
\put(37,17.4){\makebox(0,0){$\frown$}}
}
}
\put(39,35){\makebox(0,0){$H_2$}}
\put(39.5,22){\makebox(0,0){$\mathcal{Q}$}}
\end{picture}
\end{center}
\vspace{-0.7cm}

 & 

\vspace{1.1cm} 

 $\bullet$ Introduces a unique infinitesimal stratum $\mathcal{Q}$.
 
 $\bullet$ Invariant strictly drops.
 
 \\
 
 \hline
\vspace{0.01cm}

{\bf CASE D1)} 
\vspace{-0.3cm}
\begin{center}
\setlength{\unitlength}{.9mm}
\begin{picture}(40,40)
\allinethickness{1.4pt}
\put(0,20){\textcolor{green}{\line(1,0){40}}}
\put(20,0){\textcolor{darkpurple}{\line(0,1){40}}}
{\color{darkpurple}
\put(20,20){\circle{3.2}}
\put(20,20){\circle{2.5}}
\put(20,20){\circle{1.9}}
\put(20,20){\circle{1}}
\put(20,20){\circle{1.5}}
\put(20,20){\circle{0.4}}
\put(20,22){\makebox(0,0){$\smile$}}
\put(20,22){\makebox(0,0){$\smile$}}
\put(20,22){\makebox(0,0){$\smile$}}

\put(20,17.4){\makebox(0,0){$\frown$}}
\put(20,17.4){\makebox(0,0){$\frown$}}
\put(20,17.4){\makebox(0,0){$\frown$}}
}
\put(24,38){\makebox(0,0){$H_1$}}
\put(37,17.2){\makebox(0,0){$H_2$}}
\put(23.5,22.2){\makebox(0,0){$\mathcal{Q}$}}
\end{picture}
\end{center}
\vspace{-0.7cm}

&

\vspace{0.1cm}

\begin{center}
\setlength{\unitlength}{.9mm}
\begin{picture}(80,40)
\allinethickness{1.4pt}
\put(13,0){\textcolor{green}{\line(0,1){40}}}

\put(43,0){\textcolor{darkpurple}{\line(0,1){40}}}
\put(0,20){\textcolor{darkpurple}{\line(1,0){50}}}
{\color{lightblue}
\put(13,20){\circle{3.2}}
\put(13,20){\circle{2.5}}
\put(13,20){\circle{1.9}}
\put(13,20){\circle{1}}
\put(13,20){\circle{1.5}}
\put(13,20){\circle{0.4}}
{\color{darkpurple}
\put(15.4,20){\makebox(0,0){$($}}
\put(15.4,20){\makebox(0,0){$($}}
\put(15.4,20){\makebox(0,0){$($}}

\put(10.7,20){\makebox(0,0){$)$}}
\put(10.7,20){\makebox(0,0){$)$}}
\put(10.7,20){\makebox(0,0){$)$}}

\put(43,20.8){\makebox(0,0){$\smile$}}
\put(43,20.8){\makebox(0,0){$\smile$}}
\put(43,20.8){\makebox(0,0){$\smile$}}

\put(43,18.7){\makebox(0,0){$\frown$}}
\put(43,18.7){\makebox(0,0){$\frown$}}
\put(43,18.7){\makebox(0,0){$\frown$}}
}
}
\put(45,38){\makebox(0,0){$H_1$}}
\put(15.3,38){\makebox(0,0){$H_2$}}
\put(14.8,24.2){\makebox(0,0){$\mathcal{Q}'$}}
\put(26.75,14.5){\makebox(0,0){$\mathbb{A}^1=\mathbb{P}^1\setminus\{\mathcal{Q}'\}$}}

\end{picture}
\end{center}
\vspace{-0.7cm}
& 

\vspace{.8cm}

$\bullet$ Adds two new stratum $\mathcal{Q}'$ and $\mathbb{A}^1=\mathbb{P}^1\setminus\{\mathcal{Q}'\}$.

$\bullet$ Invariants strictly drop at any exceptional point.

\\

\hline

\vspace{0.01cm}
{\bf CASE D2)}
\vspace{-0.3cm}
\begin{center}
\setlength{\unitlength}{.9mm}
\begin{picture}(40,40)
\allinethickness{1.4pt}
\put(0,20){\textcolor{green}{\line(1,0){40}}}
\put(20,0){\textcolor{darkpurple}{\line(0,1){40}}}
{\color{green}
\put(20,20){\circle{3.2}}
\put(20,20){\circle{2.5}}
\put(20,20){\circle{1.9}}
\put(20,20){\circle{1}}
\put(20,20){\circle{1.5}}
\put(20,20){\circle{0.4}}
{\color{darkpurple}
\put(20,22.2){\makebox(0,0){$\smile$}}
\put(20,22.2){\makebox(0,0){$\smile$}}
\put(20,22.2){\makebox(0,0){$\smile$}}

\put(20,17.2){\makebox(0,0){$\frown$}}
\put(20,17.2){\makebox(0,0){$\frown$}}
\put(20,17.2){\makebox(0,0){$\frown$}}
}
}
\put(23.5,38){\makebox(0,0){$H_1$}}
\put(36.5,17.2){\makebox(0,0){$H_2$}}
\put(23.2,22){\makebox(0,0){$\mathcal{Q}$}}
\end{picture}
\end{center}
\vspace{-0.6cm}

& \vspace{0.1cm}

\begin{center}
\setlength{\unitlength}{.9mm}
\begin{picture}(50,40)
\allinethickness{1.4pt}
\put(13,0){\textcolor{green}{\line(0,1){40}}}

\put(37,0){\textcolor{darkpurple}{\line(0,1){40}}}
\put(0,20){\textcolor{green}{\line(1,0){50}}}
{\color{lightblue}
\put(37,20){\circle{3.2}}
\put(37,20){\circle{2.5}}
\put(37,20){\circle{1.9}}
\put(37,20){\circle{1}}
\put(37,20){\circle{1.5}}
\put(37,20){\circle{0.4}}
{\color{darkpurple}
\put(37,22.1){\makebox(0,0){$\smile$}}
\put(37,22.1){\makebox(0,0){$\smile$}}
\put(37,22.1){\makebox(0,0){$\smile$}}

\put(37,17.3){\makebox(0,0){$\frown$}}
\put(37,17.3){\makebox(0,0){$\frown$}}
\put(37,17.3){\makebox(0,0){$\frown$}}
}
}
\put(39.5,38){\makebox(0,0){$H_1$}}
\put(15.3,38){\makebox(0,0){$H_2$}}
\put(39.5,22.4){\makebox(0,0){$\mathcal{Q}'$}}
\end{picture}
\end{center}
\vspace{-0.6cm}

&

\vspace{0.1cm}

$\bullet$ Only adds a new infinitesimal stratum $\mathcal{Q}'$.

$\bullet$ The invariant remains equal, $\inv(\Q')=\inv(\Q)$.

$\bullet$ Leads to resolution or case D1) after finitely many quadratic transformations 
\\
\hline
\end{tabular}



\end{document}